\documentclass{amsart}

\usepackage{graphicx} 
\usepackage{macro}
\usepackage{tikz,tkz-euclide}
\usepackage{mathtools}
\usepackage{thmtools}
\usepackage{thm-restate}

\usepackage[utf8]{inputenc}
\usepackage[english]{babel}
\usepackage{biblatex}
\title{Exit path categories induced by group actions}
\author{Patrick Mayeda}
\date{November 11th, 2025}

\begin{document}

\begin{abstract}
    We prove a structural result concerning the exit path category associated to a manifold $M$ equipped with a smooth action of a finite group $G$.  Specifically, the functor $\Pi: \exit(M) \rightarrow \exit(M/G)$ is a right fibration and $\mathsf{Enter}(M/G)$ is classified by a natural functor $\mathsf{Enter}(M/G) \rightarrow \cO_G$, where $\cO_G$ is the orbit category of $G$.  The main technical result manipulates exit paths to immediately exiting paths, enabling lifts of homotopies in $M/G$ to homotopies in $M$.
\end{abstract}
\maketitle
\section*{Introduction}
Let $G$ be a finite group acting freely on a Hausdorff space $M$.  The quotient map $\pi: M \rightarrow M/G$ is a covering space.  This covering space is classified by a functor $\pi_1(M/G) \rightarrow BG$ insofar as there is a canonical pullback square among groupoids
\[\begin{tikzcd}
	{\pi_1(M)} && {EG} \\
	\\
	{\pi_1(M/G)} && {BG}
	\arrow[from=1-1, to=1-3]
	\arrow[from=1-1, to=3-1]
	\arrow["\lrcorner"{anchor=center, pos=0.125}, draw=none, from=1-1, to=3-3]
	\arrow[from=1-3, to=3-3]
	\arrow[from=3-1, to=3-3]
\end{tikzcd}\]
where $BG$ is the one object groupoid with $\hom(*, *) \cong G$.  A natural question arises.
\begin{q}
    If the action $G \lacts M$ is \emph{not} free, does there exist a similar classification result?
\end{q}
This work answers this question in the affirmative, at the cost of assuming that $M$ is a manifold on which $G$ acts smoothly.  The relaxation of the freeness of the action introduces nontrivial fixed points.  Pointedly, the quotient map no longer enjoys the property of being a covering space.  To address this issue, we introduce \bit{stratifications} of $M$ and $M/G$ induced by the group action.  For $x \in M$, the natural assignment of $x$ to the conjugacy class of its stabilizer subgroup in the poset of conjugacy classes is continuous.  Moreover, this descends to a continuous map from the quotient $M/G$.  These continuous maps allow the constructions of the \bit{exit path categories}, $\exit(M)$ and $\exit(M/G)$ whose objects are the points of $M, M/G$, and whose morphisms are the \bit{exit paths} up to \bit{exit path homotopy}.  An exit path is a path that respects the structure of the stratification, an exit path homotopy is a homotopy through exit paths.  $\exit(M/G)^{\op}$ is naturally given the name $\mathsf{Enter}(M/G)$.  
\par
It is through the enter/exit path category framework that a classification theorem is established.  We must also replace $BG$ with the \bit{orbit category of $G$}.  This is denoted $\cO_G$ and is the full subcategory of $\gset$ consisting of the $G$-sets with transitive group actions.  We replace $EG$ with the pointed orbit category $\cO_{G, *}$.  This is the full subcategory of pointed $\gset$ consisting of the pointed $G$-sets with transitive group actions.  There is an evident forgetful functor between the two.
\begin{theorem}
    Let $G$ be a finite group acting smoothly on a manifold $M$.  The functor $\Pi^{\op}: \mathsf{Enter}(M) \rightarrow \mathsf{Enter}(M/G)$ is recovered as a pullback of the forgetful functor $\cO_{G, *} \rightarrow \cO_G$ along a classifying functor $\mathsf{Enter}(M/G) \rightarrow \cO_G$.
    \[\begin{tikzcd}
    	{\mathsf{Enter}(M)} && {\cO_{G, *}} \\
    	\\
    	{\mathsf{Enter}(M/G)} && {\cO_G}
    	\arrow[from=1-1, to=1-3]
    	\arrow["{\Pi^{\op}}"', from=1-1, to=3-1]
    	\arrow["\lrcorner"{anchor=center, pos=0.125}, draw=none, from=1-1, to=3-3]
    	\arrow["{\mathsf{forget}}", from=1-3, to=3-3]
    	\arrow[from=3-1, to=3-3]
    \end{tikzcd}\]
\end{theorem}
\par
One hopes to recover the pullback square of groupoids from this theorem, in fact one must do so.  As such, in the case of a free action $G$, there should be an equality $\mathsf{Enter}(M) = \pi_1(M), \mathsf{Enter}(M/G) = \pi_1(M/G)$.  This is so: for a free action, for all $x \in M$, the stabilizer subgroup is the trivial group and the stratification is trivial.  The restriction of morphisms to exit paths is no longer relevant and the enter path category has the same data as the fundamental groupoid.  As for $BG$, there is a natural inclusion of $BG$ into $\cO_{G, *}$ where the unique object $*$ in $BG$ is mapped to the group $G$ itself, with the free action given by multiplication on the left.  The classifying functor $\mathsf{Enter}(M/G) \rightarrow \cO_G$ will factor through $BG$ if $G$ acts freely on $M$.
\par
To establish the classifying functor $\mathsf{Enter}(M/G) \rightarrow \cO_G$, one is in need of the following lemma.
\begin{lemma}
    \emph{(see Theorem 5.1.4)} The functor $\Pi: \exit(M) \rightarrow \exit(M/G)$ is a \bit{right fibration}
\end{lemma}
This puts us in a position to take advantage of an equivalence of categories between right fibrations over a category $C$ and presheaves on $C$.  The functor corresponding to $\Pi$ under this equivalence has as its codomain the category of sets, but one then notices that it factors through the orbit category.  The bulk of the paper is dedicated to establishing this property for the functor $\Pi$.  Specifically, given some exit path $f$ in $M/G$ with $f(0) = [x]$ and $f(1) = [y]$, fix some choice $y \in M$ as a lift for $[y]$.  First one must show that this determines a lift $\tilde{f}$ of $f$.  This is a fact that holds true outside the current context of exit paths.  Generally, for $G$ acting discontinuously on a Hausdorff space $X$, the quotient map $q: X \rightarrow X/G$ has the path lifting property, a proof of which may be found in Brown \cite{brown}.  The difficulty comes in showing that such lifts are well defined.  If $f, g: I \rightarrow M/G$ are exit path homotopic in the quotient $M/G$, the lifts $\tilde{f}, \tilde{g}$ must also be exit path homotopic in $M$.  
\par
Homotopies in general can be very unruly.  By unruly, we mean that the induced stratification of the square by post-composition with the stratifying map of $X$ can be poorly behaved, even if one limits themselves to exit path homotopies.  This prompts the introduction of \bit{immediately exiting} paths.  These are exit paths $f$ which take image in exactly two strata, $f(0) = p$ and $f(t) = q$ for all $t \in (0, 1]$ with $p < q$.  It also prompts the idea of \bit{immediately exiting homotopies}, which are homotopies of immediately exiting paths.  The stratifications on $I \times I$ induced by such homotopies are very simple.  One stratum is dedicated to $\{0\} \times I$, and another to $(0, 1] \times I$.
\begin{lemma}
    \emph{(See Theorem 4.3.10)} Let $G$ be a finite group acting smoothly on a manifold $M$.
    \begin{enumerate}
        \item Every exit path $f$ in $M/G$ is homotopic an immediately exiting path $f'$ in $M/G$
        \item Every exit path homotopy $\cH$ in $M/G$ is homotopic to an immediately exiting homotopy $\tilde{\cH}$ in $M/G$.
        \item Immediately exiting homotopies lift. The homotopy witnessing $f$ as homotopic to immediately exiting $f'$ in $M/G$ lifts.
    \end{enumerate}
\end{lemma}
This is the main technical result of the paper.  The main idea is quite simple.  One has an exiting homotopy $\cH$ relating two exit paths $f$ and $g$.  Construct a homotopy $\Phi$ from $\cH$ to $\tilde{\cH}$ where $\tilde{\cH}$ is immediately exiting between paths $f'$ and $g'$.  In particular, $\Phi$ witnesses homotopies $\cF$ between $f$ and $f'$ and $\cG$ between $g$ and $g'$.  Lifting $\cH$ is unreasonable, but one may lift $\cF, \cG$ and $\tilde{\cH}$ and rely on transitivity.
\par
The restriction to smooth actions of finite groups on manifolds means that, locally, the action looks like a linear action; the induced stratification by the poset of conjugacy classes is then \bit{conically smooth} in the sense of Ayala, Francis and Tanaka \cite{ayala}. This works to endow the strata $X_p$ of the space $X$ with collar neighborhoods $U_p$ that deformation retract onto $X_p$.  Furthermore, the maps $U_p \rightarrow X_p$, $U_p \setminus X_p \rightarrow X_p$ assemble into fiber bundles.  By way of these fiber bundles, we peel a homotopy $\cH$ off of strata step by step.  The fact that $\Pi$ is a right fibration follows readily.  
\par
We now provide a general outline of the paper.  Section $\S1$ starts off with a definition of orbit categories and provides explicit descriptions of their structure.  In addition, it elucidates some properties of orbit categories that make computation of examples for finite groups $G$ relatively easy. Section $\S2$ defines the exit path category and shows that this is indeed a category.  After doing so, we discuss the specific stratifications at hand, those induced by a finite group action on a Hausdorff space $X$.  Throughout the first half of the paper, the assumption that the space be a manifold is unnecessary.  Section $\S3$ provides a proof of lifts of exit paths.  This is an alternate proof to the one found in Brown that is specific to exit paths and makes use of the following fact that we also prove: when restricted to a single stratum associated to a conjugacy class of stabilizers, the quotient map $\pi: X \rightarrow X/G$ is a covering space.  
\par
Moving on to section $\S4$, we begin to assume that the space we are working with is a manifold and that the action is smooth.  We define immediately exiting paths and immediately exiting homotopies before going on to discuss the nature of the fiber bundles $U_p \rightarrow X_p, U_p \setminus X_p \rightarrow X_p$ previously alluded to, which arise as the fiberwise cones of a fiber bundle over a stratum $X_p$ with total space given by the \bit{link of the stratum}, denoted $\Link_{X}(X_p)$.  The discussion of fiberwise cones of fiber bundles illuminates the deformation retract $U_p \times I \rightarrow U_p$ and is an essential tool in ensuring that we may lift the exit the path homotopies $\cF, \cG$ in Lemma 0.0.3.  The process of constructing the homotopy $\Phi$ between $\cH$ and immediately exiting $\tilde{\cH}$ is given inductively by first producing a homotopy $\Phi'$ from $\cH$ to $\cH'$ where $\cH'$ takes image in a single stratum.  This is done stratum by stratum, using the fiber bundle $U_p \setminus X_p$ to peel the homotopy off the lowest dimensional stratum.  The homotopy $\Phi$ is produced by deforming a face of the homotopy cube $\Phi'$ back to the starting point of the paths in $\cH$.  This deformation is ensured by making sure the the starting face of the homotopy cube stays within some open set.  A conically smooth stratification is in particular \bit{conical}, and it is this structure which permits this deformation.  Finally, in section $\S5$, we establish $\Pi$ as a right fibration and prove the main theorem.

\section*{Relations to other works}
This work draws inspiration from the several previous works.  \emph{Systems of Fixed Point Sets} by Elmendorf \cite{elmendorf} identified the role of orbit categories in equivariant topology.  It is from this that one suspects the generalization of $BG$ to be the orbit category.  We establish the path lifting property for quotients of Hausdorff spaces in Section $\S3$.  As stated, this utilizes similar techniques to the techniques found in \cite{brown}.  The terms exit path categories were first introduced by MacPherson and defined in \cite{treumann}.  The technical process of lifting homotopies in section $\S4$ makes use of immediately exiting paths and immediately exiting homotopies, which is inspired by the $\infty$-category of definition A.6.2 in  \cite{lurie}.  The technical challenge of Lemma 0.0.4 can be rephrased by saying that the exit path category is equivalent to the immediately exiting path category.  In section $\S 2$, we show that the former can be naively taken as a category, the latter naturally takes the structure of an $\infty$-category.  (Immediately exiting) Exit path $\infty$-categories have seen considerable development such as in \cite{ayala}, \cite{barwick2020exodromy}, \cite{haine2024exodromyconicality}.  Most pertinently, \cite{ayala} establishes the existence of links between strata of stratified spaces.  This appears as Theorem 4.2.1 and plays an important role in the arguments to come.

\section*{Acknowledgements}
This article is the outcome of a project supervised by David Ayala while he was supported by the National Science Foundation under award 1945639. I am indebted to the aforementioned Doctor Ayala for all the enlightening conversations. A thank you to Jane and Katie for tolerating me.  Another thank you to the Georgia contingent for all their insight.


\section{Orbit Categories}
\subsection*{Definition, properties of orbit categories}
Every transitive group action $G \curvearrowright X$, is isomorphic as an object in the category $\gset$ to $G/H$ for some subgroup $H$.
\begin{definition}
    For a group $G$ the \emph{\textbf{orbit category of G}} denoted $\cO_G$ has as its objects $G/H$ for every $H \subseteq G$, and as its morphisms $G$-equivariant set maps.
\end{definition}
\begin{remark}
    This is a variation of the definition of the orbit category as in the introduction, where one takes $\cO_G$ to be the full subcategory of $\gset$ consisting of the transitive group actions.  This description of $\cO_G$ with explicit objects $G/H$ is sometimes written as $\mathsf{Orb}_G$.  In this notation, there is an inclusion $\orbit_G \hookrightarrow \cO_G$ that is manifestly fully faithful and essentially surjective.  Given this equivalence of categories, we give the two the same notation.
\end{remark}

\begin{prop}
    For some map $\phi: G/H \rightarrow X$ of $G$-sets, $\phi$ is completely determined by $\phi(eH)$.
\end{prop}
\begin{proof}
    $\phi(gH) = \phi(geH) = g \cdot \phi(eH)$
\end{proof}

\begin{prop}
    There is an ismorphism (of sets) $\hom_{\gset}(G/H, X) \cong X^H := \{x \in X : \forall h \in H, hx = x \}$.
\end{prop}
\begin{proof}
    Let $\phi: G/H \rightarrow X$ be some such map.  It is determined completely by the image of $eH \in X$.  One hopes that $\phi(eH) \in X^H$, the image of the coset $eH$ in $X$ is fixed under the action of $G$ for all $h \in H$.  This follows from the equation
    \begin{equation}
        \phi(eH) = \phi(h \cdot eH) = h \cdot \phi(eH)
    \end{equation}
    This establishes a set map in one direction:
    \[\begin{tikzcd}
	{\hom_{\gset}(G/H, X)} &&& {X^H} \\
	\phi &&& {\phi(eH)}
	\arrow[from=1-1, to=1-4]
	\arrow[from=2-1, to=2-4]
    \end{tikzcd}\]
    To establish the bijection, construct an inverse.  To specify the data of a $G$ equivariant map from $G/H$ to $X$, it is sufficient to choose an image for $eH$ under a map $\phi$.  The equalities in equation 1.0.1 however tell us that we may not pick \emph{any} element of $X$, it must be invariant under the action of elements of $H$.  This establishes the claim.
\end{proof}

\begin{definition}
    An \bit{endomorphism-isomorphism} category, \bit{EI-category} for short, is a category where every endomorphism is an isomorphism.
\end{definition}
\begin{prop}
    For objects $c, d$ in an EI category $C$, if there exists an isomorphism $f: c \rightarrow d$, then $\hom_C(c, d)$ consists of all isomorphisms.
\end{prop}
\begin{proof}
    Let $C$ be an EI-category and objects $c, d$ objects and $f: c \rightarrow d$ an isomorphism with inverse $f^{-1}: d \rightarrow c$.  Let  $g: c \rightarrow d$ be arbitrary.  Then $g \circ f^{-1}$ is an endomorphism of $d$, as such, it is equipped with an inverse, call it $h$.  The equation $hgf^{-1} = 1_d$ holds. Apply $f^{-1}$ on the left and $f$ on the right, then $f^{-1}hgf^{-1}f = f^{-1}f = 1_c$, so it must be so that $(f^{-1}h)g = 1_c$.  This suggests $f^{-1}h: d \rightarrow c$ as a candidate for the inverse to $g$.  It has been shown to be a left inverse, to see that it is a right inverse: $g(f^{-1}h) = (gf^{-1})h = 1_d$
    by associativity and since $h$ is an inverse to $gf^{-1}$.
\end{proof}

\begin{prop}
    $\cO_G$ is an $EI$-category. Every set $hom_{\cO_G}(G/H, G/H)$ carries a group structure isomorphic to the $W_G(H)$, the Weil group of $H$.
    $$
    W_G(H) := N_G(H)/H \cong \hom(G/H, G/H)
    $$
\end{prop}
\begin{proof}
    Suppose that $\phi: G/H \rightarrow G/H$ is an endomorphism of the object $G/H$ in $\cO_G$.  Suppose that $\phi(eH) = gH$.  Consider then the endomorphism generated by $\psi(eH) = g^{-1}H$. This defines a $G$-equivariant map if and only if $g^{-1}H$ is fixed under the action of $G$ by elements of $H$.  We have the supposition that $gH$ is fixed by these elements.  $\forall h \in H$, $hgH = gH$ holds if and only if $g^{-1}hg \in H$.  The inverse of this element, is $gh^{-1}g^{-1}$.  Since every element of $H$ is expressible as an inverse, it must be that for all $h \in H$, $hg^{-1}H = g^{-1}H$, as needed.  Computing the composition of these maps yields the morphism that sends $eH$ to itself.  Since morphisms are uniquely determined by the image of $eH$, this morphism must be the identity and these are inverses.
    \par
    $\hom(G/H, G/H)$ has been shown to be a group.  To see that it is the Weil group, construct a group homomorphism from $N_G(H) = \{g \in G : g^{-1}Hg = H \} \rightarrow \hom(G/H, G/H)$ with kernel $H$.  Since $g^{-1}Hg = H$ if and only if $g \in (G/H)^H$, there is a surjective assignment 
    \[\begin{tikzcd}
    	{N_G(H)} &&& {\hom(G/H, G/H)} \\
    	g &&& \begin{array}{c} \phi: G/H \rightarrow G/H \text{ s.t.}\\ \phi(eH) = gH \end{array}
    	\arrow["\Phi", from=1-1, to=1-4]
    	\arrow[from=2-1, to=2-4]
    \end{tikzcd}\]
    This is a group homorphism with kernel $H$, yielding the desired isomorphism.
\end{proof}

\begin{prop}
    There exists an isomorphism in $\hom_{\cO_G}(G/H, G/K)$ if and only if $H$ is conjugate to $K$.
\end{prop}
\begin{proof}
    Suppose that $\phi: G/H \rightarrow G/K$ is an isomorphism with $\phi(eH) = gK$.  Its inverse is necessarily characterized by $\phi^{-1}(eH) = g^{-1}H$.  Since these define $G$-equivariant morphisms, $\forall h \in H$, $hgK = gK$ and for all $k \in K$, $kg^{-1}H = g^{-1}H$.  The first implies that $g^{-1}Hg \subseteq K$, and the second implies that $gKg^{-1} \subseteq H$, yielding $K, H$ as conjugate.
    \par
    If $H$ is conjugate to $K$, it is enough to find a single isomorphism from $G/H$ to $G/K$.  For some $g \in G$, $g^{-1}Hg = K, g^{-1}Kg = H$.  But it is so that $kgH = gH$ if and only if $g^{-1}kg \in H$, and $hg^{-1}K = g^{-1}K$ if and only if $ghg^{-1} \in K$, these hold and we may define $\phi: G/H \rightarrow G/K$ by sending $eH$ to $g^{-1}K$ and its inverse from $G/K$ to $G/H$ by sending $eK$ to $gH$.
\end{proof}
\begin{cor}
    $\hom_{\cO_G}(G/H, G/K) \neq \emptyset$ if and only if $H$ is conjugate to some subgroup of $K$.
\end{cor}

\subsection*{Examples of orbit categories, transitive $\gset$ and pointed transitive $\gset$}
The above properties equip us with the tools to see what some of these things are for some small groups.  Consider abelian groups first.  For any abelian group $G$, no two distinct subgroups are conjugate to each other.  In this context corollary 1.0.8 becomes $\hom_{\cO_G}(G/H, G/K) \neq 0$ if and only if $H \subseteq K$.  Generally, for $G$ abelian, the structure of the category will be the subgroup lattice with a few extra arrows.  It is also quite reasonable to describe exactly what these arrows are.  All subgroups are normal in an abelian group.  If $H$ is normal in $G$, then the Weil group $W_G(H) = G/H$.  The group of endomorphisms on each object will then simply be the object itself!
\par
Furthermore, if $H \subseteq K$ are two subgroups of an abelian group, then for \emph{any} coset $gK$, $gK \in (G/K)^H$: $hgK = ghK = gK$.  Every coset is then a valid choice for the image of $eH$ for a prospective map $\phi: G/H \rightarrow G/K$.  Therefore, $|\hom_{\cO_G}(G/H, G/K)| = [G:K]$, the number of cosets in $G/K$.
\par
These observations turn the prospect of writing down an orbit category associated to an abelian group into something very reasonable. As an example, the Klein four group $K = C_2 \oplus C_2 = \{(0, 0), (1, 0), (0, 1), (1, 1)\}$, with addition component-wise mod 2.  There are three non-trivial subgroups each cyclic of order two, $C_{10}, C_{01}, C_{11}$, corresponding to the subgroup generated by $(i, j)$.  $\orbit_K$ is given below.
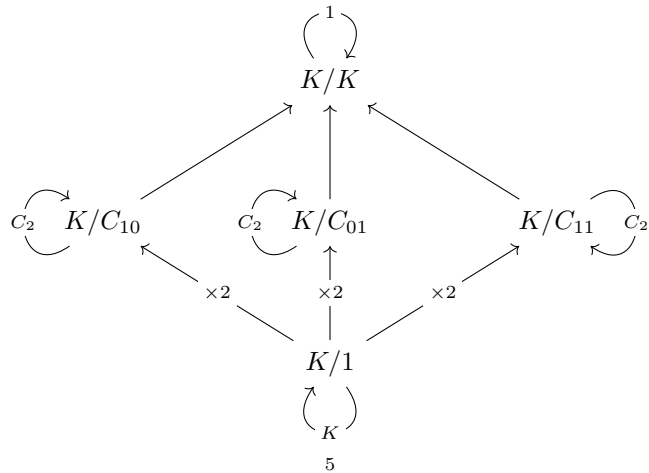
\[\begin{tikzcd}
	&& {K/K} \\
	\\
	{K/C_{10}} && {K/C_{01}} && {K/C_{11}} \\
	\\
	&& {K/1}
	\arrow["1"{description}, from=1-3, to=1-3, loop, in=55, out=125, distance=10mm]
	\arrow[from=3-1, to=1-3]
	\arrow["{C_2}"{description}, from=3-1, to=3-1, loop, in=145, out=215, distance=10mm]
	\arrow[from=3-3, to=1-3]
	\arrow["{C_2}"{description}, from=3-3, to=3-3, loop, in=145, out=215, distance=10mm]
	\arrow[from=3-5, to=1-3]
	\arrow["{C_2}"{description}, from=3-5, to=3-5, loop, in=325, out=35, distance=10mm]
	\arrow["\times2"{description}, from=5-3, to=3-1]
	\arrow["\times2"{description}, from=5-3, to=3-3]
	\arrow["\times2"{description}, from=5-3, to=3-5]
	\arrow["K"{description}, from=5-3, to=5-3, loop, in=235, out=305, distance=10mm]
\end{tikzcd}\]
This is a specific example of an abelian group, but the ideas hold in complete generality.  If one squints and sees all the hom sets as only having a singleton, it is the subgroup lattice.  After unsquinting, one can fill in the number of morphisms in each hom set without too much strain.
\par
As an example of something non-abelian.  Consider $\Sigma_3$, the symmetric group on three elements.  The subgroup lattice remains a good enough start.  There are four non-trivial subgroups, one of order three: $\langle(123)\rangle$, and three of order 2: $\langle(12)\rangle$, $\langle(13)\rangle$, $\langle(23)\rangle$.  All the subgroups of order 2 are conjugate.  Note that the size of the hom-sets $\hom_{\cO_G}(\Sigma_3/1, -)$ and $\hom_{\cO_G}(- , \Sigma_3/\Sigma_3)$ have no mystery to them.  In the first case, everything is fixed by the action of the identity.  The cardinality of the hom-set will then be the cardinality of the codomain.  In the second case $\Sigma_3/\Sigma_3$ is the singleton set, so there is only one choice of map.  A computation reveals that there is a unique map between each of the order two subgroups.  The Weil subgroup for each order two subgroup is the trivial group we leave these out of the diagram.  $\orbit_{\Sigma_3}$ is given by 
\[\begin{tikzcd}
	&& {\Sigma_3/\Sigma_3} \\
	\\
	{\Sigma_3/\langle(123)\rangle} &&& {\Sigma_3/\langle(12)\rangle} & {\Sigma_3/\langle(13)\rangle} & {\Sigma_3/\langle(23)\rangle} \\
	\\
	&& {\Sigma_3/1}
	\arrow["1"{description}, from=1-3, to=1-3, loop, in=50, out=130, distance=15mm]
	\arrow[from=3-1, to=1-3]
	\arrow["{C_3}"{description}, from=3-1, to=3-1, loop, in=140, out=220, distance=15mm]
	\arrow[from=3-4, to=1-3]
	\arrow[tail reversed, from=3-4, to=3-5]
	\arrow[from=3-5, to=1-3]
	\arrow[tail reversed, from=3-5, to=3-6]
	\arrow[from=3-6, to=1-3]
	\arrow["\times2"{description}, from=5-3, to=3-1]
	\arrow["\times3"{description}, from=5-3, to=3-4]
	\arrow["\times3"{description}, from=5-3, to=3-5]
	\arrow["\times3"{description}, from=5-3, to=3-6]
	\arrow["{\Sigma_3}"{description}, from=5-3, to=5-3, loop, in=230, out=310, distance=15mm]
\end{tikzcd}\]

\begin{definition}
    Denote the category of all transitive $\gset$ as $\cO_G$, that is the full subcategory of $\gset$, $G \curvearrowright X$ where $G$ acts transitively on $X$.  A \emph{\textbf{pointed $G$-set}} is a set with a group action and a distinguished element. The category of pointed $G$-sets has as its morphisms $G$-equivariant set maps which preserve the base point and is denoted $\gset_*$.  Denote the full subcategory of $\gset_*$ consisting of those pointed $G$-sets with transitive actions as $\cO_{G, *}$.  From here on, by the orbit category $\cO_G$ and the pointed orbit category $\cO_{G, *}$, we will mean these definitions.
\end{definition}

\begin{prop}
The evident forgetful functors assemble into a commutative diagram:
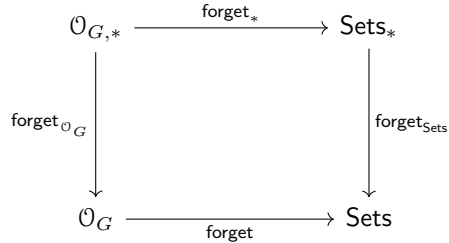
\[\begin{tikzcd}
	{\mathcal{O}_{G, *}} &&& {\mathsf{Sets}_*} \\
	\\
	\\
	{\mathcal{O}_G} &&& {\mathsf{Sets}}
	\arrow["{\mathsf{forget}_*}", from=1-1, to=1-4]
	\arrow["{\mathsf{forget}_{\cO_G}}"', from=1-1, to=4-1]
	\arrow["{\mathsf{forget}_{\mathsf{Sets}}}", from=1-4, to=4-4]
	\arrow["{\mathsf{forget}}"', from=4-1, to=4-4]
\end{tikzcd}\]
This square is a pullback of categories.
\end{prop}
\begin{proof}
    To show this, we compute the pullback of $\mathsf{forget}$ along $\mathsf{forget}_{\set}$ and notice that it is an exact description of pointed transitive $G$-sets.  An object of the pullback is a pair of objects $(X, (A, a))$ from $\cO_G$ and $\set_*$ respectively such that $\mathsf{forget}(X) = \mathsf{forget}_{\set}(A, a)$.  This is so if and only if $X = A$.  It follows that we may describe an object of the pullback as $(X, x)$.  A morphism is a pair of morphisms, $\phi: X \rightarrow Y$, a $G$-equivariant morphism in $\cO_G$, and $g: (X, x) \rightarrow (Y, y)$ a morphism of pointed sets such that the underlying set map is the same.  This is exactly a $G$-equivariant set map between $X, Y$ that preserves the base point, a morphism in $\cO_{G, *}$
\end{proof}


\section{Stratifications and Exit Path Categories}
\subsection{Properties of exit paths}
\begin{definition}
    Given a poset $P$,the \emph{\textbf{Alexandrov topology}} on $P$ is the topology consisting of the upward closed sets.  $U \subseteq P$ is open if and only if $p \in U$ and $p \leq q$ implies that $q \in U$.  Note that for some fixed $p \in P$, $\uparrow p = \{q \in P : p \leq q \}$ forms a basis for this topology.
\end{definition}
\begin{definition}
    A \emph{\textbf{stratification}} of a topological space $X$ is a continuous map $\cS: X \rightarrow P$ where $P$ is some poset with the Alexandrov topology.
\end{definition}
\begin{definition}
    A \emph{\textbf{stratum}} of a stratification $\cS: X \rightarrow P$ is the preimage of a point $p \in P$, $\cS^{-1}(p) \subseteq X$.  It is denoted as $X_p$
\end{definition}

This is the topic of the rest of the paper.  For our purposes, we take all posets to be finite and $X$ to be Hausdorff.  The rest of the section is dedicated to some useful facts about exit path categories.

\begin{definition}
    An \textbf{\emph{exit path}} in a stratified space $X$ is a path $f: I \rightarrow X$ such that for each $s \leq t$, $\cS(f(s)) \leq \cS(f(t))$. 
\end{definition}
\begin{definition}
    An \textbf{\emph{exit path homotopy}} between exit paths $f, g$ is a homotopy rel endpoints $\cH: I \times I$ such that for each $s \in I$, $\cH|_{I \times \{s\}}$ is an exit path.  Two exit paths $f, g$ are exit path homotopic (rel endpoints), written $f \simeq_{\text{exit}} g$ if there exists an exit path homotopy between the two.  
\end{definition}
\begin{definition}
    The \textbf{\emph{exit path category}} of a stratified topological space $S: X \rightarrow P$, denoted $\exit(X)$ is defined by
    \par
    i) Objects: points of $X$
    \par
    ii) Morphisms: $\hom(x, y)$ is then given by exit paths from $x$ to $y$ under the equivalence relation of exit path homotopy.  Composition is given by concatenation of paths.  The identity is the constant path.
\end{definition}
It is required (though it is perhaps 'obvious') to show that these data satisfy the axioms of a category.  To do this, it is sufficient to show that the familiar properties of homotopy hold when we restrict to exit path homotopies.
\begin{prop}
    If $f, g$ are exit paths with $f(1) = g(0)$, then $f * g$ defined by
    \begin{equation}
    \begin{cases}
        f(2t) & \text{ if } t \in [0, \frac{1}{2}] \\
        g(2t-1) & \text{ if } t \in [\frac12, 1]
    \end{cases}
    \end{equation}
    is an exit path.
\end{prop}
\begin{proof}
    Let $t_0 < t_1$. Do this in cases.  If $t_0, t_1 \in [0, \frac{1}{2}]$, then $\cS(f*g(t_i)) = \cS(f(2t_i))$.  Since $f$ is exiting and $t_0 < t_1$ implies $2t_0 < 2t_1$, the relation $\cS(f*g(t_0)) < \cS(f*g(t_1))$ holds.  The same is true for $t_0, t_1 \in [\frac{1}{2}, 1]$ since $t_0 < t_1$ also implies that $2t_0-1 < 2t_1 - 1$.  In the last case, Since $f(1) = g(0)$, for all $t_0 \in [0, \frac{1}{2}], t_1 \in [\frac12, 1]$ it is so that $\cS(f(2t_0)) \leq \cS(f(1)) = \cS(g(0)) \leq \cS(g(2t_1-1))$.
\end{proof}
The following proposition is a restatement of Theorem 3.1 in Rotman \cite{rotman} applied to exit paths and the proof goes through in the same way with one caveat, one needs to check that the proposed homotopy between $f_0*g_0$ and $f_1*g_1$ is an exit path homotopy.
\begin{prop}
    Suppose that $f_0, f_1, g_0, g_1$ are exit paths in $X$ with $f_0 \simeq_{\text{exit}} f_1$ and $g_0 \simeq_{\text{exit}} g_1$ rel endpoints.  If $f_0(1) = g_0(1), f_1(0) = g_1(0)$, then $f_0 * g_0 \simeq_{\text{exit}} f_1 * g_1$.
\end{prop}
\begin{proof}
    Let $\cF: f_0 \simeq_{\text{exit}} f_1, \cG: g_0 \simeq_{\text{exit}} g_1$ be the exit path homotopies relating each pair of paths.  The homotopy in Rotman relating the concatenations without regard to exit paths is
    \[
    \cH(t, s) := 
    \begin{cases}
        \cF(2t, s) & \text{ if }  t \in [0, \frac{1}{2}] \\
        \cG(2t-1, s) & \text{ if } t \in [\frac{1}{2}, 1] \\
    \end{cases}
    \]
    The restriction to a fixed horizontal strip, $\cH|_{I \times \{s\}}$ is realized as the concatenation $\cF|_{I \times \{s\}} * \cG|_{I \times \{s\}}$.  Since both of these are assumed to be exit paths, the concatenation is by proposition 2.1.6 and $\cH$ therefore assembles into an exit path homotopy.
\end{proof}
This proposition means that composition in our would be category of exit paths is well defined.  To finish the discussion, we'd like to show that the constant path serves as an identity and that composition is associative.  This is the statement of Theorem 3.2 in Rotman, parts i) and ii) in the exit path framework.
\begin{prop}
    For $x \in X$ with $\cS: X \rightarrow P$ a stratified space, denote the constant path at $x$ as $c_x$.  If $f(0) = x$ and $f(1) = y$, then $c_x * f \simeq_{\text{exit}} f \simeq_{\text{exit}} f * c_y$.  Associativity holds whenever possible.
\end{prop}
\begin{proof}
    Both of these follow from an examination of the homotopies witnessing these statements in the non-exit path case.  In either case, one reparametrizes the interval in a way that is both linear and monotonic across the homotopy.
\end{proof}
\begin{cor}
    The definition of exit path category as stated in 2.1.5 is in fact a category.
\end{cor}
\begin{remark}
    It is worthwile to note that part iii) of Theorem 3.2 in Rotman does \emph{not} hold.  The reverse path is not an inverse, because generally the reverse path is not an exit path.  The only case where the reverse path is admitted as an exit path is if $f: I \rightarrow X$ has its image contained in a single stratum.  In this light, one can view this framework as a generalization of the fundamental groupoid of $X$.  Taking $P = \{x\}$, then \emph{all} paths are exit paths and \emph{all} homotopies are exit path homotopies thereby yielding an equality of $\exit(X)$ with the fundamental groupoid of $X$.
\end{remark}

\begin{definition}
    The \emph{\textbf{length}} of an exit path is the number of strata it traverses, explicitly for $f: I \rightarrow X$ an exit path $\text{length}(f) := |\cS(f(I))|$.
\end{definition}

For some exit path $f: I \rightarrow X$, post-compose with the stratification of $X$.  This is a stratification of the interval.  Moreover, the image of $\cS \circ f$ is a subposet $P'$ of the poset $P$ that is linearly ordered.  This is forced by the exiting condition and the fact the the interval is linearly ordered.  If length$(f) = n$, then $P'$ is of the form $p_1 < p_2 < \dots < p_n$, where $\cS(f(0)) = p_1$ and $\cS(f(1)) = p_n$.

\begin{prop}
    Let $P = p_1 < \dots < p_n$ be a linearly ordered poset of size $n$.  For $\cS: I \rightarrow P$ a surjective stratification such that $s \leq t$ implies $\cS(s) \leq \cS(t)$, then $I$ has $n$ strata of the form $[a_0 = 0, a_1], (a_1, a_2], \dots, (a_{n-1}, a_n = 1]$.
\end{prop}
\begin{proof}
    By induction on the length of the poset.  If the poset is a singleton, there is necessarily 1 stratum of the form $[0, 1]$.  Suppose the above holds for linearly ordered posets up to lengh $n$ and suppose that $P$ is of length $n+1$, $P = p_0 < p_1 < \dots < p_n$ and consider the preimage of the closed set $\downarrow p_{n-1}$.  This is a closed set in the unit interval and we claim it to be of the form $[0, a]$ for some $a \in I$.  
    \par
    Define $a := \sup{\{t \in I\} : \cS(t) \in \downarrow p_{n-1}\}}$.  Since $\cS(0) = p_0$, this set is nonempty, and so $a$ exists.  It is manifestly a limit point of the set $\cS^{-1}(\downarrow p_{n-1})$ and is therefore an element thereof by closedness.  Every element of $\cS^{-1}(\downarrow p_{n-1})$ is in $[0, a]$ since $a$ is the supremum of this set and $0$ is necessarily the infimum.  Suppose then that $x \in (0, a)$ and $\cS(x) = p_n$.  This yields $x < a$ and $\cS(a) = p_{n-1} < p_n = \cS(x)$, contradicting monotonicity.  
    \[\begin{tikzcd}
    	I && {[0, a]} && {\downarrow p_{n-1}} \\
    	t && at
    	\arrow["\phi", from=1-1, to=1-3]
    	\arrow["{\mathcal{S}|}", from=1-3, to=1-5]
    	\arrow[from=2-1, to=2-3]
    \end{tikzcd}\]
    is now a surjective stratification of $I$ by a linearly ordered poset of length $n$.  It therefore has the form $[0, b_1], \dots, (b_{n-1}, 1]$.  This is to say that $\phi^{-1}(\cS^{-1}(p_k)) = (b_k, b_{k+1}]$ (excepting the case of $p_0$ where the left endpoint is closed, though it proceeds in the same way). $\phi$ a homeomorphism implies $\cS^{-1}(p_k) = \phi((b_k, b_{k+1}]) = (\phi(b_k), \phi(b_{k+1})]$
    With $\cS^{-1}(p_n) = I \setminus \cS^{-1}(\downarrow p_{n-1}) = (a, 1]$, this gives the desired structure of the strata.
\end{proof}

\begin{definition}
    Call the stratification of $I$ induced by an exit path as in proposition 2.1.7 a \emph{\textbf{segmentation}} of $I$ and call the $[a_k, a_{k+1}]$ (with both endpoints closed) the \emph{\textbf{segments}} of $f$.
\end{definition}

This is a convenience.  For $f$ an exiting path in some quotient of a space, we'll want to lift $f$ to an exit path in the original space.  With this, it will suffice to lift one segment at a time and ensure that they agree on the points where they overlap.

\subsection{Stratifications by group actions}

\begin{prop}
    Suppose that $G$ is a group acting on a Hausdorff topological space $X$.  If $S_G$ is the subgroup lattice, then $\cS: X \rightarrow S_G^{\op}$ by $\cS(x) = G_x$ the stabilizer of $x$ is a stratification of $X$.
\end{prop}
\begin{proof}
    There is a basis of closed sets for the topology on $S_G^{\op}$ given by the downsets on each $H \in S_G^{\op}$ which we denote $\downarrow H := \{K \subseteq G : K \supseteq H\}$.  The op induces the reversal of the inclusion.  It is enough to check that $f^{-1}(\downarrow H)$ is closed for each $H$.
    $$
    f^{-1}(\downarrow H) = \{x \in X : f(x) = G_x \supseteq H \} = X^H
    $$
    where $X^H := \{x \in X : \forall h \in H, \; h \cdot x = x\}$.  Indeed, if the stabilizer of $x$ has $H$ as a subgroup, then certainly all the elements of $H$ fix $x$ under the group action.
    \par
    The problem now reduces to showing that for each $H \subseteq G$, $X^H$ is closed in $X$.  There is a homeomorphism $\psi: \Delta \rightarrow X$ where $\Delta := \{(x, x)\} \subseteq X \times X$ is the diagonal.    Furthermore, continuity of the group action implies that for each $h \in H$ the map $x \rightarrow h\cdot x$ is continuous which in turn means that the graph of each map $\Gamma_h = \{(x, h\cdot x)\} \subseteq X \times X$ is closed in the product.  The intersection $\Gamma_h \cap \Delta$ is then the set of points that are fixed under the action of $h$.  To find the set of points which are fixed under the action of $h$ for \emph{every} $h \in H$, it is then necessary to take the intersection over all such sets.
    $$
    F := \{(x, x) : \forall h \in H, h \cdot x = x \} = \cap_{h \in H} (\Gamma_h \cap \Delta) \subseteq X \times X
    $$
    As the intersection of a number of closed sets, $F$ is closed in $X \times X$.  $F \subseteq \Delta$ then yields that then that $F$ is in fact closed in $\Delta$.  But it is evident that $\psi(F) = X^H$ meaning that $X^H$ is closed.
\end{proof}
Note here that we use the fact that $X$ is Hausdorff when we take each $\Gamma_h$ to be closed in the product, but we do \emph{not} use the fact that $G$ is finite.  Realizing $X^H$ as a closed set depends on an \emph{intersection} of closed sets indexed by the number of elements of $H$.  Of course, $H$ need not be finite for this intersection to remain closed.
\par
The ultimate goal is a stratification of $X$ \emph{and} $X/G$, so this should factor through the quotient.  It doesn't, so we pass to $P_G^{\op}$, the poset of conjugacy classes of subgroups.
\[\begin{tikzcd}
	X && {S_G^{\op}} && {P_G^{\op}} \\
	\\
	& {X/G}
	\arrow[from=1-1, to=1-3]
	\arrow["\pi"', from=1-1, to=3-2]
	\arrow[from=1-3, to=1-5]
	\arrow["{\exists! \tilde{\cS}}"', dashed, from=3-2, to=1-5]
\end{tikzcd}\]
That the continuous map $\tilde{\cS}$ in the above diagram exists, is given by the following: 
\begin{prop}
    $x \in \text{Orbit}(y)$ implies that $G_x$ is conjugate to $G_y$.
\end{prop}
\begin{proof}
    Let $x \in \text{Orbit}(y)$ and $g \in G$ such that $x = g \cdot y$ and $g^{-1} \cdot x = y$.  Let $k \in G_x$.
    \[
    (g^{-1}kg) \cdot y = g^{-1}\cdot (k \cdot( g \cdot y)) = g^{-1}\cdot (k \cdot x) = g^{-1} \cdot x = y
    \]
    This implies that $g^{-1}G_x g \subseteq G_y$.
    \par
    For the reverse inclusion, let $h \in G_y$.
    \[
    y = h \cdot y = (hg^{-1}) \cdot x
    \]
    Applying the action of $g$ to both sides of the equation yields $x = (ghg^{-1}) \cdot x$.  For all $h \in G_y$, it is so that $ghg^{-1} \subseteq G_x$.  This yields $G_y \subseteq g^{-1}G_xg$ and subsequently the equality $g^{-1}G_xg = G_y$.
\end{proof}
\begin{cor}
    For $G \curvearrowright X$ a continuous group action on a Hausdorff topological space $X$, both $X$ and $X/G$ have natural stratifications by the poset $P_G^{\op}$ of conjugacy classes of subgroups of $G$ given by 
    \[\begin{tikzcd}
    	x &&&& {G_x} \\
    	X &&&& {P_G^{\op}} \\
    	&&&& {[G_x]} \\
    	{X/G} \\
    	{[x]}
    	\arrow[from=1-1, to=1-5]
    	\arrow["{\mathcal{S}}", from=2-1, to=2-5]
    	\arrow["\pi", from=2-1, to=4-1]
    	\arrow["{\tilde{\mathcal{S}}}"', from=4-1, to=2-5]
    	\arrow[from=5-1, to=3-5]
    \end{tikzcd}\]
\end{cor}

Fixing an arbitrary topological space $P$, we can consider the category over $P$ with objects $X \rightarrow P$ and morphisms $\phi: X \rightarrow Y$ such that the triangle
\[\begin{tikzcd}
	X &&& P \\
	\\
	Y
	\arrow[from=1-1, to=1-4]
	\arrow["\phi"', from=1-1, to=3-1]
	\arrow[from=3-1, to=1-4]
\end{tikzcd}\]
commutes.  When $P$ is a poset, these objects are of course stratifications.  One can hope that passing to the enter path category is functorial.  This is the content of the next proposition.

\begin{prop}
    In the case that $P$ is a poset, the assignment on objects $(X \rightarrow P) \rightarrow \exit(X)$ assmebles into a functor $\exit(-): \Top/P \rightarrow \Cat$.
\end{prop}
\begin{proof}
    The assignment on objects is given.  For some $\phi: X \rightarrow Y$ with the appropriate commutativity conditions, we are in need of a functor between the two categories $\exit(\phi): \exit(X) \rightarrow \exit(Y)$.  The objects of each of these categories are the points of $x$ and the points of $Y$.  The only sensible choice of the assigment on objects for the functor $\exit(\phi)$ is $x \rightarrow \phi(x)$.  Given some exit path $f$ starting at $x_0$ and ending at $x_1$, commutativity in the diagram
\[\begin{tikzcd}
	I &&& X &&& P \\
	\\
	&&& Y
	\arrow["f", from=1-1, to=1-4]
	\arrow["{\phi \circ f}"', dashed, from=1-1, to=3-4]
	\arrow["{\mathcal{S}_X}", from=1-4, to=1-7]
	\arrow["\phi"', from=1-4, to=3-4]
	\arrow["{\mathcal{S}_Y}"', from=3-4, to=1-7]
\end{tikzcd}\]
    gives for all $t$ the equality $\cS_Y(\phi(f(t))) = \cS_X(f(t))$.  The composition $\phi \circ f$ is then necessarily an exit path starting at $\phi(x_0)$ and ending at $\phi(x_1)$.  This is the assignment on morphisms.  
    \par
    We are considering exit paths up to homotopy.  To for this to be well defined, if $\cH$ witnesses an exit path homotopy from $f$ to $f'$, then there must be an exit path homopy $\cH'$ from $\phi \circ f$ to $\phi \circ f'$. Modify the diagram, replacing $f$ with our homotopy
    \[\begin{tikzcd}
    	{I \times I} &&& X &&& P \\
    	\\
    	&&& Y
    	\arrow["{\mathcal{H}}", from=1-1, to=1-4]
    	\arrow["{\phi \circ \mathcal{H}}"', dashed, from=1-1, to=3-4]
    	\arrow["{\mathcal{S}_X}", from=1-4, to=1-7]
    	\arrow["\phi"', from=1-4, to=3-4]
    	\arrow["{\mathcal{S}_Y}"', from=3-4, to=1-7]
    \end{tikzcd}\]
    $\phi \circ \cH: I \times I \rightarrow Y$ is the desired homotopy.  This is by the eqality $(\phi \circ \cH)|_{I \times \{\alpha\}} = \phi \circ (\cH|_{I \times \{\alpha\}})$
\end{proof}

\section{Lifting exit paths}
\subsection{Restrictions to stratum are covering spaces}
First, some facts about group actions of finite groups on Hausdorff spaces are needed.
\begin{definition}
    A group $G$ acts on a space $X$ \emph{\textbf{discontinuously}} if, for all $x \in X$, there exists some open set $V$ about $x$ such that, for all $g \notin G_x$, $g \cdot V \cap V = \emptyset$
\end{definition}
\begin{prop}
    If $G$ is a finite group acting on a Hausdorff space $X$, $G \curvearrowright X$, it necessarily acts discontinuously.
\end{prop}
\begin{proof}
    Let $x \in X$.  Enumerate the elements of $g \notin G_x$ as $g_1, \dots g_n$ such that each $g_i \cdot x$ is distinct.  This yields $n+1$ distinct points in $X$, $x, g_1\cdot x, \dots, g_n \cdot x$.  Since $X$ is Hausdorff, about each of these points, there exists disjoint open sets $V_0, V_1, \dots V_n$ separating the points, i.e. such that $V_i \cap V_j = \emptyset$ for $i \neq j$.  Multiplication by $g_i$ is a homeomorphism, so $g_1 \cdot V_0$ is an open set about $g_1 \cdot x$, as is $U_1 := (g_1 \cdot V_0) \cap V_1$.  Doing this for all $V_i$ and taking $U_0 = V_0$, we have constructed $U_i$ disjoint about each $x_i$ such that $U_i \subseteq g_i \cdot V_i$.  For each $U_i$, $g^{-1}_i \cdot U_i$ is an open set about $x_0$.  Define $U := \cap g^{-1}_i \cdot U_i$.  Since $g_i \cdot U \subseteq U_i \subseteq V_i$, this yields the desired open set about $x_0$ witnessing the fact that $G$ acts on $X$ discontinuously.
\end{proof}
\begin{prop}
    If $G$ acts continuously on $X$, then the quotient map $\pi: X \rightarrow X/G$ is an open map.
\end{prop}
\begin{proof}
    Let $U$ be an open subset of $X$.  It suffices to show that $\pi^{-1}(\pi(U))$ is open in $X$.  This is the set of all elements $x \in X$ such that there exists $y \in U$ and $g \in G$ with $g \cdot y = X$.  This is the claim
    $$
    \pi^{-1}(\pi(U)) = \cup_{g \in G} g \cdot U
    $$
    Each multiplication $g \cdot -$ is a homeomorphism, and so each $g \cdot U$ is an open set.  This is then open as union of open sets.
\end{proof}

These are general facts, the specific lemmas needed for our endeavor come down to guaranteeing the existence of certain well behaved open sets $U$ of $X$.  Ultimately, for some arbitrary $x \in X$ we'd like to know that given \emph{any} $U$ open about $x$, $\exists V \subseteq U$ with some particular properties.  We give names to these properties and put them in a single proposition.
\begin{lemma}
Let $G \curvearrowright X$ be a continuous group action of finite $G$ on Hausdorff $X$.  For any open set $U$ about a point $x \in X$, $\exists V \subseteq U$ such that all of the following hold:
    \begin{enumerate}
        \item (Inclusion) $\forall y \in V, G_y \subseteq G_x$.
        \item (Discontinuity) $\forall g \notin G_x$, $g \cdot V \cap V = \emptyset$.
        \item (Symmetry) $\forall g \in G_x$, $g \cdot V = V$.
    \end{enumerate} 
Furthermore, for any open set $V$ about a point $x$ satisfying conditions 2 and 3, the following is true: $\pi^{-1}(\pi(V)) = \coprod g_k \cdot V$, where the $g_k$ in the cofactors of this disjoint union range over a single element from each of the cosets of $G/G_x$.
\end{lemma}
\begin{proof}
If, for each open $U$ about some point $x$, there exists $V_1, V_2, V_3 \subseteq U$ satisfying inclusion, discontinuity and symmetry respectively, we may take the intersection $V_1 \cap V_2 \cap V_3 \subseteq U$ and this will satisfy all three conditions simultaneously.  It then suffices to show the the existence of such a $V$ one property at a time.
\par
\emph{Existence of (1):} The assignment $x \rightarrow G_x$ is always a continuous map between $X$ and $S_G^{\op}$, as in the proof of proposition 2.2.1.  $\uparrow G_x := \{H : H \subseteq G_x\}$ is an open set of $S_G^{\op}$.  Take the intersection of $U$ with the preimage of this open set as $V$.
\par
\emph{Existence of (2):} By proposition 3.1.2, $G$ acts on $X$ discontinuously $\implies \exists W$ about $x$ such that $\forall g \notin G_x$, $g \cdot W \cap W = \emptyset$.  Take $U \cap W$.
\par
\emph{Existence of (3):} For $U$ some open neighborhood of $x$, consider $\cap_{g \in G_x} g \cdot U$.  Note that since $G_x$ is finite, this is an open set.  Furthermore, $x$ is in this intersection, so it is a non-empty neighborhood of $x$.  Letting $h \in G_x$, we claim that $h \cdot (\cap_{g \in G_x} g \cdot U) = (\cap_{g \in G_x} g \cdot U)$.  To see this, that the map $h \cdot -: X \rightarrow X$ is a homeomorphism, specifically it is injective $\implies h \cdot (\cap_{g \in G_x} g \cdot U) = \cap_{g \in G_x} hg \cdot U$.  Now, since $h$ is in fact an element of $G_x$, this is just the same as reindexing this intersection: $\cap_{g \in G_x} hg \cdot U = \cap_{k \in G_x} k \cdot U$.
\par
For the last statement, note that for any $V$ satisfying the symmetry condition, if $h, k$ are elements of the same coset of $G/G_x, hG_x = kG_x \iff h^{-1}k \in G_x$, then $h^{-1}k \cdot V = V \iff k \cdot V = h \cdot V$.  That is, the image of $V$ under the maps $h \cdot -$ and $k \cdot -$ are the same.  It is generally so that $\pi^{-1}(\pi(V)) = \cup_{g \in G} g \cdot V$, this observation that two elements of the same coset produce the same set in the union, allows us to toss out all but one $g \cdot V$ from each coset.  The discontinuity condition guarantees that this isn't just a union, but a disjoint union, justifying the use of the coproduct symbol instead of union.
\end{proof}

With the information of proposition 2.1.13, an exit path $f: I \rightarrow X/G$ pulls back to a stratification of $I$ of the form $[0, a_1], (a_1, a_2], \dots, (a_{n-1}, 1]$.  Each half open segment $(a_{k-1}, a_k]$ takes image in a single stratum.  Expressing strata as covering spaces would, given a lift of $f(a_k)$, allow us to lift all of the half open interval.  One hopes then to extend the lift uniquely as the limit point in $X$ of our lift.  This lifts the entire segment $[a_{k-1}, a_k]$ and moreover determines a lift of $f(a_{k-1})$ to be used to lift the next segment. This is the plan that we will execute.  First we establish some notation.  For $G \curvearrowright X$ our continuous group action as usual, $X$ and $X/G$ are stratified by functions in the commutative diagram
\[\begin{tikzcd}
	X &&& {P_G^{\op}} \\
	\\
	& {X/G}
	\arrow["{\mathcal{S}}", from=1-1, to=1-4]
	\arrow["\pi"', from=1-1, to=3-2]
	\arrow["{\tilde{\mathcal{S}}}"', from=3-2, to=1-4]
\end{tikzcd}\]
For $[K]$ some conjugacy class of subgroups of $P$, denote the stratum $\cS^{-1}([K])$ in $X$ by $X_K$ and the stratum $\cS^{-1}([K])$ in $X/G$ by $X/G_K$. Commutativity of the diagram implies that $\pi(X_K) = X/G_K$ and $\pi^{-1}(X/G)_K = X_K$.
\begin{prop}
    The restriction of $\pi$ to $X_K$ and its image $X/G_K$ is a covering space with $|Orbit(x)|$ sheets, for $x$ any element of $X_K$.    
\end{prop}
\begin{proof}
    Before anything, it is unambiguous to say that we have $|\text{Orbit}(x)|$ sheets: for any $x \in X_K$, the stabilizer of $x$ is conjugate to $K$.  Therefore $|\text{Orbit}(x)| = |G/G_x| = |G/K|$.
    \par
    For $x \in X_K$, take some open subset $V \subseteq X$ of $x$ satisfying the properties in lemma 3.1.4. (this starts with some open $U$ and finds $V$ in $U$, so if necessary just start with all of $X$).  $U := V \cap X_K$ is then an open neighborhood of $x \in X_K$ with some other nice features.  The inclusion condition strengthens to an equality.  Indeed, $\forall y \in X_K$, $G_y$ conjugate is to $G_x$.  Then $\forall y \in U$, $G_y \subseteq G_x$ and conjugacy implies equality.  It is also worth noting that multiplication by $g \in G$ defines an action on $X_K$.  $\forall x \in X_K$, $g \cdot x \in Orbit(x)$ has a stabilizer that is conjugate to $K$, it is therefore true that the image of the set $X_K$ under each homeomorphism $g \cdot -$ is $X_K$ itself.
    \par
    For convenience, give the restriction of $\pi$ its own name: $p: X_K \rightarrow X/G_K$.  Then $p|_U: U \rightarrow p(U)$ is a bijection.  If $y \neq z$ are such that $p(y) = p(z)$, this is so if and only if $y, z$ are in the same orbit, that is $\exists g \in G$ such that $g \cdot y = z$.  This is to say that $g \notin G_y = G_z = G_x$, but this supplies us with $z \in g \cdot V \cap V$, a contradiction.
    \par
    Better still, $p|_U$ is an open map, and therefore a homeomorphism.  By proposition 3.1.3 $\pi$ is an open map.  By the observation made earlier, $p$ is the restriction of $\pi$ to a saturated subspace.  Any such restriction of an open map remains an open map.  $p|_U$ is the further restriction of an open map to an open set.  This is also necessarily another open map.
    \par
    The properties 1,2, and 3 descend to the open set $U \subseteq X_K$.  Inclusion and discontinuity are immediate.  Symmetry follows from the fact that for all $y \in U$, $G_y = G_x$, so for any $g \in G_x$, $g \cdot U = U$ in the strictest sense possible, nothing has moved.  It follows then that $p^{-1}(p(U))$, the set consisting of all those elements in $X_K$ in the orbit of some $y \in U$ is exactly $\coprod g \cdot U$ taken over a choice of element $g$ from each coset.  $g \cdot -$ is a homeomorphism, $g \cdot U \cong U \cong p(U)$ and all of these disjoint finishes the proof.
\end{proof}


\subsection{Lifting exit paths across strata} The proof of the fact that $p: X_K \rightarrow X/G_K$ is a covering space is to use the following theorem found in Rotman \cite{rotman}:
\begin{theorem}
    (Theorem 10.13, Rotman) Let $Y$ be connected and locally path connected, and let $f: (Y, y_0) \rightarrow (X, x_0)$ be continuous.  If $(\tilde{X}, p)$ is a covering space of $X$, then there exists unique $\tilde{f}: (Y, y_0) \rightarrow (\tilde{X}, \tilde{x}_0)$ (where $\tilde{x}_0 \in p^{-1}(x_0)$) lifting $f$ if and only if $f_*\pi_1(Y, y_0) \subset p_*\pi_1(\tilde{X}, \tilde{x}_0)$
\end{theorem}
The condition $f_*\pi_1(Y, y_0) \subset p_*\pi_1(\tilde{X}, \tilde{x}_0)$ is guaranteed for $Y$ with trivial $\pi_1$, as is the case with $I$ which is contractible.

\begin{definition}
    For $\cS:X \rightarrow P$ a stratified space, an \emph{\textbf{immediately exiting segment}} is, for $[a, b] \subseteq I$, a continuous map $f: [a, b] \rightarrow X$ such that $\tilde{S}(f(a)) = p$ and, for all $t \in (a, b]$, $\tilde{S}(f(t)) = q$ with $p < q$
\end{definition}
\begin{remark}
    Note that by proposition 2.1.7 any exit path $f: I \rightarrow X/G$ of length $n$ splits $I$ up into $n-1$ immediately exiting segments and one segment that stays in a single stratum.  Our goal is to lift such an exit path to $\tilde{f}: I \rightarrow X$, as such, it suffices to lift each one of these immediately exiting segments and glue together along the points where they overlap.
\end{remark}
\begin{prop}
    (Lifts of immediately exiting segments) For $f: [a, b] \rightarrow X/G$ an immediately exiting segment, if $f(b) = [x]$ and $x_1$ is a choice of lift for $[x]$, then $\exists! \tilde{f}: I \rightarrow X$ such that $f(b) = x_1$ and $\pi \circ \tilde{f} = f$.
\end{prop}
\begin{proof}
    Since $f|_{(a, b]}$ has image in a single stratum, the stated theorem means that a choice of a lift for $f(b)$ determines a lift $\tilde{f}$ of this restriction.
    \par
    We'd like now to extend $\tilde{f}$ to $\tilde{f}_{ext}: [a, b] \rightarrow X$ in such a way so that $\tilde{f}_{ext}$ defines a lift for $f$ on the entire interval.  To start, if $f(a) = [y]$, enumerate the lifts of $[y]$ as $y_1, \dots, y_n$.  For $y_1$, fix an open set $V$ about $y_1$ as in lemma 3.1.4.  As in proposition 3.1.3, the quotient map is an open map, yielding $\pi(V)$ an open set about the point $[y]$ in $X/G$.  Since $f(a) = [y]$, $\exists \alpha \in [a, b]$ such that $\forall t \in [a, \alpha)$, we have $f(t) \in \pi(V)$.  It must be so then that $\forall t \in (a, \alpha)$, $\tilde{f}(t) \in \pi^{-1}(\pi(V)) = \coprod g_k \cdot V$, where the coproduct is taken over a selection of $g_k$ from each coset of $G_{y_1}$ and each of the $g_k \cdot V$ are disjoint.  The open interval $(a, \alpha)$ is connected, so the image of $(a, \alpha)$ under $\tilde{f}$ must lie in a single cofactor $g_k \cdot V$ for the image to remain connected.  Without loss of generality, assume that this cofactor is $e \cdot V$, the cofactor containing $y_1$.  Define $\tilde{f}_{ext}(a) = y_1$.
    \par
    The only point in $[a, b]$ at which $\tilde{f}_{ext}$ might not be continuous is where we have defined the extension at the beginning of the interval.  To check continuity of $\tilde{f}_{ext}$ it suffices then to check continuity just at the point $a$. With this in mind, let $U$ be some open set about $y_1$. With $V$ as in the previous paragraph, again invoking lemma 3.1.4, $\exists W \subseteq U \cap V$ satisfying the symmetry and discontuity conditions.  We do not know if $\tilde{f}_{ext}$ is continuous, but we do know that as set functions, there is an equality $\pi \circ \tilde{f}_{ext} = f$.  Applying the open quotient map and using continuity of $f$, there exists an open set $A \subseteq I$ about $a$ such that $f(A) \subseteq \pi(W)$.  We may take the connected component in which $a$ lies and assume that $A$ is of the form $[a, \omega)$.  Define $\gamma = \text{min}(\alpha, \omega)$. Commutativity as set functions implies that 
    $$
    \tilde{f}_{ext}([a, \gamma)) \subseteq \pi^{-1}(\pi(W)) = \coprod g_k \cdot W
    $$
    We know however that $\tilde{f}_{ext}(a) = y_1 \in W$.  Connectedness of $(a, \gamma)$ and continuity of $\tilde{f}$ implies that $\tilde{f}_{ext}((a, \gamma))$ is in one cofactor.  Furthermore, we know that $\tilde{f}((a, \gamma)) \subseteq V$.  The condition that $W \subseteq V$ implies that the containment $g_k \cdot W \subseteq g_k \cdot V$ holds for every $g_k$ and further implies that $g_k \cdot W \cap V = \emptyset$ for every $g_k$ that isn't in the stabilizer of $y_1$.  We'd like to show that $\tilde{f}_{ext}((a, \gamma)) = \tilde{f}((a, \gamma)) \subseteq W$.  For the sake of contradiction, suppose that it is contained in a different cofactor.  This means that $\tilde{f}((a, \gamma)) \subseteq g_k \cdot W \cap V = \emptyset$, clearly an impossibility.
    \par
    This shows that $\tilde{f}_{ext}([a, \gamma)) \subseteq W \subseteq U$ and witnesses continuity of the extension at the point $a$.
    \par
    Uniqueness of the lift along $(a, b]$ is guaranteed by theorem 3.3.1.  Uniqueness of the extension is guaranteed by the construction. We did so by showing that the image of $\tilde{f}$ before we extended it was in exactly one cofactor of the preimage and there was exactly one lift in each cofactor by which we could define the extension.
\end{proof}
\begin{prop}
    (Lifts of exit paths) For an exit path $f: I \rightarrow X/G$ and a choice of a lift of $f(1)$, $\exists!$ lift $\tilde{f}: I \rightarrow X$ satisfying $\pi \circ \tilde{f} = f$.  These are \emph{all} the lifts of $f$.
\end{prop}
\begin{proof}
    $I \rightarrow X/G \rightarrow P$ provides a segmentation of $I$ into the closed segments $[0, a_1], [a_1, a_2], \dots, [a_n, 1]$.  By the previous lemma, a choice of lift for $f(1)$ determines a unique lift of the segment $[a_n, 1]$.  This in particular provides a choice of lift for $f(a_n)$ which further provides a lift of the segment $[a_{n-1}, a_n]$ that crucially agrees on the end point with the lift over the segment $[a_n, 1]$.  These two lifts glue together to form a lift of all of $[a_{n-1}, 1]$.  Proceed inductively.  Upon getting to the last segment, there is no need to invoke Lemma 3.3.2, the entire segment is in a single stratum and we may appeal to theorem 3.3.1.  These lifts are all distinct as they all disagree at the endpoint of the interval by construction.
    \par
    To see that these are all the lifts of $f$, any lift $\tilde{f}$ of $f$ must have $\tilde{f}(1)$ evaluate to some lift of $f(1)$ and therefore must be the unique lift corresponding to $\tilde{f}(1)$.
\end{proof}

\section{Lifting Homotopies}
To show that $\Pi$ is a fibration, it is not enough to lift exit paths.  One needs to be assured that the lift is well defined.  Given two exit paths $f \simeq_{\text{exit}} g$ in $X/G$, one needs to be able to say that the two lifts $\tilde{f}, \tilde{g}$ are exit path homotopic.  Looking at the exit path homotopy $\cH$ between $f$ and $g$, a choice of a lift of $f(1)$ and $g(1)$ determines a choice of lift for \emph{every} horizontal strip, $\cH|_{I \times \{s\}}$.  This is clearly a lift of $\cH$ as a set function.  One can hope that this is continuous, but continuity on the horizontal strips does not guarantee continuity on the entirety of the square.  After this idea, one can resort to constructing a lift in the same way we lifted paths, stratum by stratum.  Morally, one should expect more success doing this.  This is taking the success in one dimension and applying it to the next.  The problem however is that the strata in $I \times I$ induced by the stratification of $X/G$ have very few restraints.  There is not a decomposition of $I \times I$ into half-open half-closed intervals as in the case of $I$.  One may even think of homotopies inducing fractal boundaries.  With this in mind, it is difficult to envision the success of an application of covering space theory. 
\par
At this juncture we impose an extra hypothesis, namely that the stratifications of $X, X/G$ are conically smooth.  This is met if $X$ is a manifold and the action of $G$ is smooth. The definition of a conically smooth stratification is inductive and rather involved; it may be found in Ayala, Francis and Tanaka's work on stratified spaces \cite{ayala}.  The utility of this extra hypothesis as it pertains to lifting homotopies is the existence of tubular neighborhoods about strata that deform nicely onto the strata themselves.  More than a deformation, they assemble into fiber bundles over the strata. Through these mechanisms, we construct a lift.

\begin{definition}
    An \emph{\textbf{immediately exiting path}} in a stratified space $\cS: X \rightarrow P$ is an exit path $f: I \rightarrow X$ such that $\cS(f(0)) = p$ and, for all $t \in (0, 1]$, $\cS(f(t)) = q$ with $p < 1$.
\end{definition}
\begin{definition}
    An \emph{\textbf{immediately exiting homotopy}} (rel endpoints) is an exit path homotopy $\cH: I \times I \rightarrow X$ such that, for all $s \in I$, $\cH|_{I \times \{s\}}$ is an immediately exiting path. 
\end{definition}
\begin{prop}
    Suppose that $\cH: I \times I \rightarrow X/G$ is an immediately exiting homotopy between immediately exiting paths $f, g$.  A choice of lift of $f(1) = g(1)$ uniquely determines a choice of lift $\tilde{\cH}$ of $\cH$.
\end{prop}
\begin{proof}
    It follows that when one pulls back the stratification of $X$ to $I \times I$ by way of an immediately exiting homotopy that there are two strata: $\{0\} \times I$ associated to $p$ and $(0, 1] \times I$ associated to $q$.
    \par
    Use Theorem 3.3.1 to lift $(0, 1] \times I$, the stratum associated to $q > p$.  Call this lift $\tilde{\cH}|$.  Now, $\cH(0, s) = [x] \in X/G$ for all $s \in I$.  Let $x_1, \dots, x_n$ be an enumeration of the lifts of $[x]$.  Let $U$ be a disjoint, symmetric open neighborhood about $x_1$ as in Lemma 3.1.4.  $\pi(U)$ is open about $[x]$.  For all $s \in I$ there exists an open set $V_s$ about $(0, s)$ such that $\cH(V_s) \subseteq \pi(U)$.  We may patch these into a connected open cover of $\{0\} \times I$.  Taking the union of all these $V_s$ to be $V$, we have a connected open set $V$ containing $\{0\} \times I$ with the feature that $\cH(V) \subseteq \pi(U)$.  Since $V \setminus \{0\} \times I$ remains an open and connected subset of $(0, 1] \times I$, it must be so that $\tilde{\cH}|(V \setminus \{0\} \times I)$ takes image in exactly one of the $\pi^{-1}(\pi(U)) = \coprod g \cdot U$.  Suppose without loss of generality that it takes image in $U$ and extend $\tilde{\cH}|$ to all of $I \times I$ by defining $\tilde{\cH}(0, s) := x_1$.
    \par
    This extension is continuous everywhere except perhaps on the left edge of the square.  Let $W$ be some arbitrary open set about $x_1 = \tilde{\cH}(0, s)$.  Again, find some $W_1 \subseteq W$ disjoint and symmetric about $x_1$ as in Lemma 3.1.4.  In fact, we may take $W_1 \subseteq W \cap U$ for $U$ the open set about $x_1$ in the previous paragraph. $\cH$ remains continuous at $(0, s)$, so there must be some $V_1$ about $(0, s)$ such that $\cH(V_1) \subseteq \pi(W_1)$.  We may take this to be connected.  $V_1 \setminus \{0\} \times I$ is remains connected in $(0, 1] \times I$.  Continuity of $\tilde{\cH}|$ then yields that $\tilde{\cH}(V_1 \setminus \{0\} \times I)$ is a subset of a single cofactor of $\coprod g \cdot W_1$.  Since we judiciously shrunk $W_1$ to be contained entirely inside of $U$, it follows that it must be the cofactor $W_1$.  Since all of the things in $V_1 \cap \{0\} \times I$ map to $x_1 \in W_1 \subseteq W$ by definition, we are finished.
\end{proof}
This is inspiration for the road ahead.  Though the induced stratification of $I \times I$ need not have reasonable strata, if the strata \emph{are} reasonable one has hopes of seeing a lift using a similar argument go through.  Note that the assumption that $\cH$ is an immediately exiting homotopy is too strong and that the conclusion, one may lift immediately exiting homotopies by using covering space theory, is not strong enough.  Nevertheless, we place this here as motivation for the proceedings to follow.  To begin, we give a definition of \emph{conically stratified} as in Lurie \cite{lurie}.  We go on to develop some lemmata about cones, fiberwise cones of fiber bundles and their stratifications.

\subsection{Fiberwise cones of fiber bundles}
\begin{definition}
    For a topological space $L$, the \emph{\textbf{cone on $L$}} is the topological space given by the pushout
    \[\begin{tikzcd}
    	{L \times\{0\}} && {L \times [0, 1)} \\
    	\\
    	{\{0\}} && {C(L)}
    	\arrow[hook, from=1-1, to=1-3]
    	\arrow["{!}"', from=1-1, to=3-1]
    	\arrow[from=1-3, to=3-3]
    	\arrow[from=3-1, to=3-3]
    	\arrow["\lrcorner"{anchor=center, pos=0.125, rotate=180}, draw=none, from=3-3, to=1-1]
    \end{tikzcd}\]
\end{definition}
\begin{definition}
    For a poset $P$, the \emph{\textbf{left cone on $P$}}, denoted $P^{\triangleleft}$, is the poset $P$ with a freely adjoined minimal element.  This is realizable as a pushout
    \[\begin{tikzcd}
    	{P \times 0} && {P \times (0 < 1)} \\
    	\\
    	0 && {P^{\triangleleft}}
    	\arrow[from=1-1, to=1-3]
    	\arrow[from=1-1, to=3-1]
    	\arrow[from=1-3, to=3-3]
    	\arrow[from=3-1, to=3-3]
    	\arrow["\lrcorner"{anchor=center, pos=0.125, rotate=180}, draw=none, from=3-3, to=1-1]
    \end{tikzcd}\]
\end{definition}
\begin{remark}
    If $\cS: L \rightarrow P$ is a stratification of the space $L$, then $C(L)$ is naturally stratified by the poset $P^{\triangleleft}$.  This follows immediately from the universal property of pushouts and a stratificaiton $s: [0, 1) \rightarrow 0 < 1$ given by assigning $s(0) = 0$ and $s(t) = 1$ for all $t \in (0, 1)$.  This stratification of $C(L)$ is witnessed as the induced $\tilde{\cS}$ in the diagram
    \[\begin{tikzcd}
	{L \times\{0\}} && {L \times [0, 1)} \\
	\\
	{\{0\}} && {C(L)} \\
	&& {P\times\{0\}} && {P \times (0 < 1)} \\
	\\
	&& {\{0\}} && {P^{\triangleleft}}
	\arrow[hook, from=1-1, to=1-3]
	\arrow["{!}"', from=1-1, to=3-1]
	\arrow["{\cS \times !}"{description}, from=1-1, to=4-3]
	\arrow[from=1-3, to=3-3]
	\arrow["{\cS \times s}", from=1-3, to=4-5]
	\arrow[from=3-1, to=3-3]
	\arrow["{!}", from=3-1, to=6-3]
	\arrow["\lrcorner"{anchor=center, pos=0.125, rotate=180}, draw=none, from=3-3, to=1-1]
	\arrow["{\exists!\tilde{\cS}}"{description}, dashed, from=3-3, to=6-5]
	\arrow[from=4-3, to=4-5]
	\arrow[from=4-3, to=6-3]
	\arrow[from=4-5, to=6-5]
	\arrow[from=6-3, to=6-5]
	\arrow["\lrcorner"{anchor=center, pos=0.125, rotate=180}, draw=none, from=6-5, to=4-3]
    \end{tikzcd}\]
    We denote the cone point of $C(L)$ as $*$.  If one removes the cone point, they are left with the subspace $L \times (0, 1)$.  $\tilde{\cS}(l, t)$ for $t \in (0, 1)$ is the same as $\cS(l)$.  The cone point is mapped to the newly adjoined minimal element.
\end{remark}
\begin{definition}
    The \emph{\textbf{$\epsilon$ cone of $L$}} is the cone on $L$ where the pushout is taken over $L \times [0, \epsilon)$, instead of $L \times (0, 1)$.  Denote this $C_\epsilon(L)$
\end{definition}
\begin{prop}
    For all $\epsilon \in [0, 1)$, there is an open embedding $C_\epsilon(L) \hookrightarrow C(L)$. If $L$ is compact, such open sets form a topological basis for $C(L)$ around the cone point: for all $U$ open about $*$, there exists $\epsilon > 0$ such that $C_\epsilon(L) \subseteq U$.
\end{prop}
\begin{proof}
    The natural map $L \times [0, 1) \rightarrow C(L)$ is a quotient map.  Indeed, one can define $C(L)$ just as well as the quotient of the equivalence relation $(l, 0) \sim (m, 0)$ for all $l, m \in L$ with no other non trivial relations.  To see that $C_\epsilon(L)$ is open, it then suffices to look at the preimage of this set in $L \times [0, 1)$.  Its preimage is $L \times [0, \epsilon)$, open in the product topology.
    \par
    Letting $U$ by some arbitrary open set about the cone point $*$, the preimage of $U$ in $L \times [0, 1)$ is open and contains the set $L \times \{0\}$.  About each point $(l, 0)$, there exists an open set contained in the preimage which may be taken to be of the form $U \times [0, \epsilon_l)$.  Over all $l \in L$, such points form a cover of $L \times \{0\}$.  Using the fact that $L$ is compact and taking the minimum over all $\epsilon_l$ in a finite cover yields an $\epsilon >0$ such that $L \times [0, \epsilon)$ is a subset of the preimage.  The image of this is $C_\epsilon(L)$.
\end{proof}

\begin{definition}
    A stratification $\cS: X \rightarrow P$ is \emph{\textbf{conical}} if, $\forall x \in X$, there exists a pointed topological space $(V, x_0)$, a stratified topological space $L \rightarrow P_{> \cS(x)}$ and an open embedding $\phi: V \times C(L) \rightarrow X$ such that 
    \[\begin{tikzcd}
    	{(x_0, *)} &&& x \\
    	{V \times C(L)} &&& X \\
    	{C(L)} \\
    	{P_{>\cS(x)}^{\triangleleft}} &&& P
    	\arrow[from=1-1, to=1-4]
    	\arrow[from=2-1, to=2-4]
    	\arrow[from=2-1, to=3-1]
    	\arrow["\cS", from=2-4, to=4-4]
    	\arrow[from=3-1, to=4-1]
    	\arrow[hook, from=4-1, to=4-4]
    \end{tikzcd}\]
    the square commutes.
\end{definition}
\begin{remark}
    Conically smooth stratificafications are in particularly conically stratified.  Moreover, for each $V$ coordinate in open basics $V \times C(L)$, one may take $V$ to be a vector space.  In what is perhaps an abuse of notation, we will refer to open basic neighborhoods $V \times C(L)$ about points $x$ in conically stratified neighborhoods.  If $X$ is a conically smooth stratification, one may take $V$ to be a vector space and the distinguished point $x_0$ to be the origin.
\end{remark}
\begin{definition}
    For a fiber bundle $\pi: E \rightarrow B$, the \emph{\textbf{fiberwise cone}} of $E$, written $C^{fib}(E)$ is given by the pushout
    \[\begin{tikzcd}
    	{E\times\{0\}} && {E \times [0, 1)} \\
    	\\
    	B && {C^{fib}(E)}
    	\arrow[from=1-1, to=1-3]
    	\arrow["\pi"', from=1-1, to=3-1]
    	\arrow[from=1-3, to=3-3]
    	\arrow["{\frak{loc}}"', from=3-1, to=3-3]
    	\arrow["\lrcorner"{anchor=center, pos=0.125, rotate=180}, draw=none, from=3-3, to=1-1]
    \end{tikzcd}\]
    where $\frak{loc}$ is the \emph{cone locus map}.  Since $\pi$ is surjective, one can also offer $C^{fib}(E)$ as a quotient of $E \times [0,1)$ by the equivalence relation $(e, 0) \sim (f, 0)$ if and only if $\pi(e) = \pi(f)$.
\end{definition}
\begin{remark}
    Embeddings are stable under pushout in the category of topological spaces.  The inclusion $E \times \{0\} \rightarrow E \times [0, 1)$ is an embedding, so $\frak{loc}$ is also an embedding.  Furthermore, if $\cS: E \rightarrow P$ is a stratification of $E$, it follows that $C^{fib}(E)$ is stratified by $P^{\triangleleft}$ in direct analogue with the natural stratification of the cone on $L$.  The projection map $L \times \{0\} \rightarrow \{0\}$ is a trivial example of a fiber bundle, expressing the fiberwise cone as a generalization of the cone on a space $L$.  With this stratification in mind, one can go further and say that $\frak{loc}$ is an embedding of $B$ that surjects onto the lowest dimensional stratum of $C^{fib}(E)$.  Having made this remark, it is also worth noting that the remaining propositions in the section on properties of fiber bundles remain true for cones on a compact space $L$.
\end{remark}

\begin{prop}
    There exists a deformation retract $\cD$ of $C^{fib}(E)$ onto its cone locus, $\frak{loc}(B)$.  Fixing $(x, s) \in C^{fib}(E) \setminus \frak{loc}(B)$, this deformation retract determines a path to $[x, 0]$ from $(x, s)$.  The reverse of this path is immediately exiting.
\end{prop}
\begin{proof}
    Recall that $- \times I$ preserves colimits.  This is to say that 
    \[\begin{tikzcd}
    	{E\times\{0\} \times I} && {E \times [0, 1) \times I} \\
    	\\
    	{B \times I} && {C^{fib}(E) \times I}
    	\arrow[from=1-1, to=1-3]
    	\arrow[from=1-1, to=3-1]
    	\arrow[from=1-3, to=3-3]
    	\arrow[from=3-1, to=3-3]
    \end{tikzcd}\]
    is a pushout diagram and we may use its mapping out property to define $\cD$.
    \[\begin{tikzcd}
    	&& {(x, t, r)} \\
    	{E\times\{0\} \times I} && {E \times [0, 1) \times I} && {(x, t(1-r))} \\
    	&&&& {E \times [0,1)} \\
    	{B \times I} && {C^{fib}(E) \times I} \\
    	\\
    	& B &&& {C^{fib}(E)}
    	\arrow[from=1-3, to=2-5]
    	\arrow[from=2-1, to=2-3]
    	\arrow[from=2-1, to=4-1]
    	\arrow[from=2-3, to=3-5]
    	\arrow[from=2-3, to=4-3]
    	\arrow[from=3-5, to=6-5]
    	\arrow[from=4-1, to=4-3]
    	\arrow[from=4-1, to=6-2]
    	\arrow["\exists!\cD"{description}, dashed, from=4-3, to=6-5]
    	\arrow[from=6-2, to=6-5]
    \end{tikzcd}\]
    It is clear that such a $\cD$ is induced, and that we are supplied with the desired deformation retract.  At $0$ it is the identity, at time $r = 1$ everything is mapped to the cone locus. The points in the cone locus remain fixed for all times $r$.  Fix $(x, s)$.  The deformation retract then defines a path from $(x, s)$ to $[x, 0]$ in the cone locus.  Looking at the reverse of this path, $[x, 0]$ is in the stratum associated to the cone point of $P^{\triangleleft}$, so all we need to show for this to be immediately exiting with respect to the stratification of $C^{fib}(E)$ is that for $r \neq 0$, $\cD(x, s, r)$ is in the same stratum of $C^{fib}(E)$ as $(x, s)$.  This is an immediate consequence of commutativity of the diagram inducing the stratification of $C^{fib}(E)$.
\end{proof}
\begin{observation}
    There isn't anything particularly salient about the fact that we fixed a \emph{point} $(x, s)$.  One may very well fix the image of a \emph{path} in $C^{fib}(E) \setminus \frak{loc}(B)$.  If this path $\gamma$ is contained in a single stratum, the result is an immediately exiting homotopy between the immediately exiting paths ending at $\gamma(0)$ and $\gamma(1)$.
\end{observation}

\begin{prop}
    If $\pi: E \rightarrow B$ is a fiber bundle, then $\pi \circ \mathsf{proj}: E \times X \rightarrow B$ is a fiber bundle with fibers $F \times X$.  If $\sim$ is an equivalence relation on $E$ such that $e \sim f \implies \pi(e) = \pi(f)$, and if $B$ is locally compact and Hausdorff, then the natural map $\omega: E/\sim \; \rightarrow B$ is a fiber bundle with fibers $F/\sim$.
\end{prop}
\begin{proof}
    For the product, let $b \in B$ and let $U$ be an open set about $b$ with $\phi$ witnessing a local trivialization.  $(\pi \circ \mathsf{proj})^{-1}(U) = \mathsf{proj}^{-1}(\pi^{-1}(U)) = \pi^{-1}(U) \times X$.  In the local trivialization, replace $\phi$ with $\phi \times \text{Id}$.
    \par
    For the quotient, the stipulation that $e \sim f \implies \pi(e) = \pi(f)$ yields the map $\omega$ from the quotient.  Furthermore, for a fixed fiber $F_b = \pi^{-1}(b)$, it allows for a sensible definition of $F_b/\sim$.  Indeed, if there are $e \in F_b$ and $f \in F_c$, there can only be a relation $e \sim f$ if $\pi(e) = b = c = \pi(f)$.  Each fiber then inherits an equivalence relation from the equivalence relation on $E$.  Let $b \in B$ and let $U$ be an open set about $b$ with $\phi: \pi^{-1}(U) \rightarrow U \times F$ witnessing the trivialization.  Now, we have an equality of $\pi^{-1}(U) = q^{-1}(\omega^{-1}(U))$ due to commutativity of the quotient diagram.  The goal then is to construct a filler $\psi$
    \[\begin{tikzcd}
    	{q^{-1}(\omega^{-1}(U)) = \pi^{-1}(U)} &&&&&&& {U \times F} \\
    	\\
    	& {\omega^{-1}(U)} &&&& {U \times F/\sim} \\
    	\\
    	&&& U
    	\arrow["\phi", from=1-1, to=1-8]
    	\arrow["{q|}"', from=1-1, to=3-2]
    	\arrow["{q_F}", from=1-8, to=3-6]
    	\arrow["\psi"{description}, dashed, from=3-2, to=3-6]
    	\arrow["{\omega|}"', from=3-2, to=5-4]
    	\arrow["{\mathsf{proj}}", from=3-6, to=5-4]
    \end{tikzcd}\]
    in the diagram that is also a homeomorphism.  This is induced from the fact that the restriction of $q$ to an open set in $q^{-1}(\omega^{-1}(U))$ is a quotient map.  Since we have assumed that $B$ is locally compact and Hausdorff, $q_F$, the product of the quotient map from $F$ to $F/\sim$ with the identity on $U$ is also a quotient map.  An inverse is induced by $\phi^{-1}$ and the universal property of $q_F$.
\end{proof}
\begin{cor}
    The two maps
    \begin{equation}
        \cD|_{C^{fib}(E) \setminus \frak{loc}(B) \times \{0\}}: C^{fib}(E)\setminus \frak{loc}(B) \rightarrow \frak{loc}(B)
    \end{equation}
    \begin{equation}
        \cD|_{C^{fib}(E) \times \{0\}}: C^{fib}(E) \rightarrow \frak{loc}(B)
    \end{equation}
    are fiber bundles with fibers $F \times (0, 1)$ and $C(F)$ respectively.  To differentiate these two restrictions, denote the first as $\cD_\circ$ and the second as $\cD_\bullet$ to evoke that the first removes the cone locus and that the second retains it.
\end{cor}
\begin{definition}
    The \emph{\textbf{$\epsilon$ fiberwise cone}} on $\pi: E \rightarrow B$ is the pushout over $E \times [0, \epsilon)$ instead of $[0, 1)$.  It is denoted $C^{fib}_\epsilon(E)$.  Denote the corresponding fiber bundles over $\frak{loc}(B)$ with total space $C^{fib}_\epsilon$ by $\cD_{\circ, \epsilon}$, $\cD_{\bullet, \epsilon}$.
\end{definition}
\begin{remark}
    One now is in pursuit of an analog of Proposition 4.1.5 in the context of fiberwise cones and not just cones.  It is evident that for all $\epsilon \in (0, 1)$ there is an open embedding of $C_\epsilon^{fib}(E)$ into $C^{fib}(E)$.  It is also evident that $C_\epsilon^{fib}(E) \rightarrow \frak{loc}(B)$ is a fiber bundle with fibers $C_\epsilon(F)$ and that $C_\epsilon^{fib}(E) \setminus \frak{loc}(B) \rightarrow \frak{loc}(B)$ is a fiber bundle with fibers $F \times (0, \epsilon)$.
\end{remark}
\begin{lemma}
    For any compact subset $C$ of $\frak{loc}(B)$ and for any open $V \subseteq C^{fib}(E)$ about $C$, there exists an $\epsilon > 0$ such that for all $c \in C$, the fibers over $c$ in $C_\epsilon^{fib}(E)$ over $c$ are contained in $V$, $D_{\bullet, \epsilon}^{-1}(c) \subseteq V$.
\end{lemma}
\begin{proof}
    For an open set $U \subseteq B$, consider sets of the form 
    $$
    U_\epsilon := \bigg\{[e, t] : e \in \pi^{-1}(U), t \in [0, \epsilon) \bigg\} \subseteq C^{fib}(E)
    $$
    Denote the natural quotient map $q: E \times [0, 1) \rightarrow C^{fib}(E)$.  It follows that $q^{-1}(U_\epsilon) = \pi^{-1}(U) \times [0, \epsilon) \subseteq E \times [0, 1)$, meaning that for all $\epsilon \in (0, 1)$, $U_\epsilon$ is open.  Let $b \in B$ be some point and let $V$ be an open set about $b$.  We claim that there exists $U \subseteq B$ and $\epsilon > 0$ such that $U_\epsilon \subseteq V$.  Now, $q^{-1}(V)$ assembles into an open set about $q^{-1}(b) = \pi^{-1}(b) \times \{0\} \cong F \times \{0\}$, where $F$ is the (compact) fiber over $b$ in the original fiber bundle $\pi: E \rightarrow B$.  About each $(f, 0) \in F \times \{0\}$, we may put an open set of the form $U_f \times [0, \epsilon_f)$ that is contained in $q^{-1}(V)$.  This assembles into a finite cover of compact $F \times \{0\}$, so we may do this with a single $\epsilon$ by taking a minimum.
    $$
    F \times \{0\} \subseteq \bigcup \big((U_{f_k} \times [0, \epsilon)\big) = \big(\bigcup U_{f_k}\big) \times [0, \epsilon) \subseteq q^{-1}(V)
    $$
    Since $q$ is surjective, it suffices to find $U \subseteq B$ an open set about $b$ such that $\pi^{-1}(U) \subseteq \cup U_{f_k}$. This would yield $U_\epsilon$ such that $q^{-1}(U_\epsilon) \subseteq q^{-1}(V)$.  Surjectivity means that applying $q$ to both sides gives $U_\epsilon \subseteq V$.
    \par
    The problem is reduced to a question about fiber bundles: Let $W$ be some open set about the fiber of $b$.  Then there exists an open set $U$ about $b$ such that $\pi^{-1}(U) \subseteq W$.  Trivialize the fiber bundle by some open set $V$ about $b$.  This reduces the problem further to the case of a projection from a product: let $W$ be an open set about $\{b\} \times F \subseteq V \times F$, then there exists $U \subseteq V$ such that $\mathsf{proj}^{-1}(U) = U \times F \subseteq W$.  This follows from compactness of $F$ and a direct application of the tube lemma.
    \par
    Back to the problem at hand, note that $C^{fib}_\epsilon(E) = U_\epsilon$ where we take $U$ to be the entirety of $B$. About each point $c$ in the compact set $c$, place an open set of the form $U_{\epsilon_c} \subseteq V$.  This assembles into a finite cover.  The evident containment $U_{\epsilon} \subseteq U_{\delta}$ for $\epsilon < \delta$ means that we can take $\epsilon$ to be the same over all open sets $U$ upon taking a minimum.  Take $C^{fib}_\epsilon(E)$ with this value of $\epsilon$.  It is clear that for each $U$ open in $B$, the inclusion $U_\epsilon \hookrightarrow C^{fib}_\epsilon(E)$ is an open embedding.  Furthermore, for all $c \in U$, $\cD_{\bullet, \epsilon}^{-1}(c)$ is contained in $U_\epsilon$.  This observation yields the result.
\end{proof}

\subsection{A stratified homotopy lifting property for conically smooth stratifications.}
Having established the properties of fiberwise cones, we now proceed to the setting of conically smooth stratifications.  We lay out the general process explicitly.  Let $\cH$ be an exit path homotopy rel endpoints.  Let $V \times C(L)$ be some open basic about $\cH(0, s)$ in $X$.
\begin{enumerate}
    \item If $X_p$ is the lowest dimensional stratum of the image of $\cH$, lift $\cH$ to $\cH'$ where $\cH'$ does \emph{not} take image in $X_p$.
    \item Construct a homotopy $\Phi$ from $\cH$ to $\cH'$ such that for all points in the homotopy variable $\alpha$, $\Phi|_{I \times I \times \{\alpha\}}$ is an exit path homotopy and such that $\Phi|_{\{0\} \times I \times \{\alpha\}} \subseteq V \times C(L)$.
    \item Proceed inductively until one has a homotopy from $\cH$ to $\tilde{\cH}$ where $\tilde{\cH}$ takes image in a single stratum.
    \item Take the left face of the homotopy cube relating $\cH$ and $\tilde{\cH}$.  It's image under the homotopy $\Phi$ is entirely in $V \times C(L)$.  Deform it back to the original point.
\end{enumerate}
The result will be a homotopy between $\cH$ and an immediately exiting homotopy.  One must also take care to ensure that the induced homotopies $\cF, \cG$ as in Theorem 0.0.2 may be lifted. The following theorem is essential in employing this line of thinking.

\begin{theorem}
    Let $\cS: X \rightarrow P$ be a conically smooth stratified space. Fix a stratum $X_p$ in $X$.  There exists a proper fiber bundle $\Link_{X_p}(X) \rightarrow X_p$ of stratified spaces where $\Link_{X_p}(X)$ is stratified by $P_{>p}$.  There exists a stratified open embedding $\phi: C^{fib}(\Link_{X_p}(X)) \rightarrow X$ such that $\phi \circ \frak{loc}$ is equal to the inclusion of $X_p$ into $X$
    \[\begin{tikzcd}
    	& {X_p} \\
    	\\
    	{C^{fib}(\Link_{X_p}(X))} &&&& X \\
    	\\
    	{P_{>p}^{\triangleleft}} &&&& P
    	\arrow["\frak{loc}"', hook', from=1-2, to=3-1]
    	\arrow["\iota", hook, from=1-2, to=3-5]
    	\arrow["{\exists\phi}"', dashed, hook, from=3-1, to=3-5]
    	\arrow[from=3-1, to=5-1]
    	\arrow[from=3-5, to=5-5]
    	\arrow[hook, from=5-1, to=5-5]
    \end{tikzcd}\]
\end{theorem}

\begin{cor}
    For any stratum $X_p$ of $X$, there exists an open neighborhood $U_p$ of $X_p$ and a deformation retract $\cD$ of $U_p$ onto $X_p$.  Restrictions of this deformation retract yield fiber bundles: $\cD_\bullet: U_p \rightarrow X_p$, $\cD_\circ: U_p \setminus X_p \rightarrow X_p$. 
\end{cor}

We lift with respect to these fiber bundles.  That they come from the deformation retract of the fiberwise cone of the link onto its cone locus should be the focus of this, specifically as it pertains to the second half of Proposition 4.1.10.  Deforming onto the cone locus yields exit paths (and homotopies) whose strata in $I \times I$ are as reasonable as one would like.  
\par
Now we proceed to the technical details of ensuring the lift as in $(1)$ of our outline exists and may be taken to be exiting.  We cite the familiar result for fiber bundles found as Proposition 4.48 in Hatcher \cite{Hatcher} as inspiration.

\begin{prop}
    A fiber bundle $\pi: E \rightarrow B$ has the homotopoy lifting property with respect to all $CW$ pairs $(A, X)$.
\end{prop}

Fiber bundles over paracompact spaces have the lifting property with respect to \emph{all} spaces, but we will not need the strength of this statement.  There is, however, a requirement to strengthen the statement cited.  Firstly, one must put problem into a situation where the fiber bundles $\cD_\circ, \cD_\bullet$ can be utilized.  An idea is encapsulated in the following diagram:
\[\begin{tikzcd}
	A &&&& {U_p \setminus X_p} \\
	\\
	{\cH^{-1}(U_p)} && {U_p} && X
	\arrow["\cH", from=1-1, to=1-5]
	\arrow[hook, from=1-1, to=3-1]
	\arrow["{\cD_\circ}", from=1-5, to=3-5]
	\arrow[dashed, from=3-1, to=1-5]
	\arrow["\cH"', from=3-1, to=3-3]
	\arrow["{\cD_\bullet}"', from=3-3, to=3-5]
\end{tikzcd}\]
for some subspace $A$ of the open set $\cH^{-1}(U_p)$.  This open set need not be nice in any sense of the word.  It \emph{is} an open set about the closed set $\cH^{-1}(X_p)$, where $X_p$ is the lowest dimensional stratum in which $\cH$ takes image.  $\cH^{-1}(X_p)$ has a special form.
\begin{definition}
    A \emph{\textbf{leftward closed}} subspace of $I \times I$ is a set $C$ such that, for all $(t, s) \in C$, if $t' \leq t$, then $(t', s) \in C$.  A \emph{\textbf{rightward closed}} set is defined analogously with the reverse inequality: $(t, s) \in C$, $t' \geq t$ implies $(t', s) \in C$.  The \emph{\textbf{leftward (rightward) closure}} of a closed set $C \subseteq I \times I$ are given by
    \[
    C_{\text{left}} = \{(t, s) : \exists t_0 \geq t, (t_0, s) \in C \}
    \]
    \[
    C_{\text{right}} = \{(t, s) : \exists t_0 \leq t, (t_0, s) \in C \}
    \]
\end{definition}

\begin{remark}
The exiting condition on $\cH$ forces $\cH^{-1}(X_p)$ to be leftward closed in $I \times I$.  Note that the definition of leftward/rightward closure need not assume that the set in question is open or closed.  We will take the rightward/leftward closure of a set to mean in general the rightward/leftward closure of a \emph{closed} set, in light of the next proposition which requires the set to be closed as a hypothesis.
\end{remark}

\begin{prop}
    The rightward (leftward) closure of a closed subset of $I \times I$ is again closed.
\end{prop}
\begin{proof}
    We only show the case of rightward closure, the case of leftward closure being completely analogous. Let $C$ be a closed subset of $I \times I$.  Suppose that $(t_k, s_k)$ is some sequence of points in $C_{\text{right}}$ converging to $(t, s)$.  This is so if and only if $t_k \rightarrow t$ and $s_k \rightarrow s$ in each coordinate.  Each point $(t_k, s_k)$ is an element of $C_{\text{right}}$ means that there are points for each $k \in \NN$, $(t_{j_k}, s_k) \in C$ with $t_{j_k} \leq t_k$.  The $t_{j_k}$ are some sequence of points in the interval, so pass to a convergent subsequence.  After reindexing, this yields $(t_l, s_l)$ where $s_l$ is a subsequence of $s_k$ and all points are in $C$ and $t_l$ converges to some $t_0 \leq t$ and $s_l$ converges to $s$.  The rightward closed condition on $C_{\text{right}}$ implies that $(t, s)$ is an element, giving closedness.
\end{proof}
\begin{theorem}
    Let $\cH$ be an exit path homotopy in a conically smooth stratified space $\cS: X \rightarrow P$.  $\cH$ need not be rel endpoints.  Let $X_p$ be the lowest dimensional stratum in which $\cH$ takes image.
    \begin{enumerate}
        \item $\exists A \subseteq A_{\text{left}} \subseteq I \times I$ such that $\cD_\circ$ has a stratified homotopy lifting property; that is, there exists a continuous filler in the diagram
        \[\begin{tikzcd}
        	A &&&& {U_p \setminus X_p} \\
        	\\
        	{A_{\text{left}}} && {U_p} && X
        	\arrow["{\cH|}", from=1-1, to=1-5]
        	\arrow[from=1-1, to=3-1]
        	\arrow["{\cD_\circ}", from=1-5, to=3-5]
        	\arrow["{\exists\text{ exiting } \tilde{\cH}}"{description}, dashed, from=3-1, to=1-5]
        	\arrow["{\cH|}"', from=3-1, to=3-3]
        	\arrow["{\cD_\bullet}"', from=3-3, to=3-5]
        \end{tikzcd}\]
        $A$ may taken to be closed and homeomorphic to a finite coproduct of intervals.
        \item There exists a homotopy $\Phi: A_{\text{left}} \times I \rightarrow X$ between $\tilde{\cH}$ and $\cH|$ such that for all $a \in A$, $\alpha \in I$, $\Phi(a, \alpha) = \tilde{\cH}(a) = \cH|(a)$, i.e. the homotopy is rel $A$.
        \item If $\cH(I \times \{0\} \cap A_{\text{left}}) \subseteq B$, for $B$ some open set, we may construct $\tilde{\cH}, \Phi$ such that for all $\alpha \in I$, $\Phi(I \times \{0\} \cap A_{\text{left}}, \alpha) \subseteq B$.
    \end{enumerate}
\end{theorem}
This is the technical formulation of parts $(1), (2)$ at the beginning of the section and serves as the necessary tool for each inductive step alluded to in $(3)$.  It is of course implicit that $A_{\text{left}} \subseteq \cH^{-1}(U_p)$.  This will take some time to prove, and much of the theorem relies on the intelligent construction of $A$.  This is the content of the next lemma.
\begin{lemma}
    Let $C$ be a leftward closed subset of the square and let $U$ be some open set with $C \subseteq U$.  Assume that $\{1\} \times I \subseteq U^C$.  There exists a piecewise linear embedding $f: \coprod_{i \in \cI} I \rightarrow I \times I$ such that 
    \begin{enumerate}
        \item $\mathsf{im}f \subseteq U \setminus C$
        \item The indexing set $\cI$ is finite
        \item $\mathsf{im}f$ separates $C$ and $U^C$:  the removal of the image of $f$ disconnects the square into at least two connected components, with $U^C$ in one and $C$ in the other.
        \item For all $s \in I$, $\mathsf{im}f \cap I \times \{s\}$ is the empty set or a singleton.
    \end{enumerate}
\end{lemma}
\begin{proof}
    Take the rightward closure of $U^C$ as $U^C_{\text{right}}$.  Since $C$ is leftward closed and $C \subseteq U$, it is true that $C \cap U^C_{\text{right}} = \emptyset$.  The open set $U_0 := U \setminus U^C_{\text{right}}$ must also contain the closed set $C$.
    \par
    Let us show that $U_0$ is leftward closed: for the sake of contradiction, suppose that there is some $(t_0, s) \notin U_0$ with $(t, s) \in U_0$ and $t_0 < t$; there are two cases: (i) $(t_0, s) \in U$, (ii) $(t_0, s) \notin U$. In the first case, $(t_0, s) \notin U$ implies that $(t_0, s) \in U^C \subseteq U^C_{\text{right}}$.  This latter set is rightward closed, as such $(t, s) \in U^C_{\text{right}}$ producing a contradiction, since $U^C_{\text{right}} \cap U_0 = \emptyset$.  In the second case, $(t_0, s) \in U$, and $(t_0, s) \notin U_0$ means that $(t_0, s)$ must be an element of $U^C_{\text{right}}$, producing the same contradiction.
    \par
    Cover compact $C$ by finitely many open boxes $\mathsf{Box}_\epsilon(c)$ such that the closure is contained in $U_0$.  This produces yet another smaller open set containing $C$ that is contained in $U_0$, 
    $$
    C \subseteq \bigcup_k \mathsf{Box}_{\epsilon_k}(c_k) \subseteq U_0 \subseteq U
    $$
    Denote the closure of the union of these boxes as $\beta$. There are inclusions $C \subseteq \beta \subset U_0$.  Take the leftward closure of this set as $\beta_{\text{left}}$.  Since $U_0$ is leftward closed, there continuous to be a containment $C \subseteq \beta \subseteq \beta_{\text{left}} \subseteq U_0$.  Furthermore, $\partial \beta_{\text{left}} \cap C = \emptyset$ and the boundary of $\beta_{\text{left}}$ separates $C$ from $U^C$.  Since $\beta$ is the union of a finite number of closed boxes, we can think of $\beta_{\text{left}}$ as a union of a finite number of leftward closed rectangles, the leftward closure of each closed box in $\beta$.  Considering now the form that $\partial \beta_{\text{left}}$ takes, it is necessarily a collection of a finite number of horizontal and vertical lines.  Let $l, m$ be two such vertical lines and let $\pi_s$ denote projection onto the $s$ coordinate of $I \times I$.  Suppose further that $\pi_s(l) \cap \pi_s(m) \neq \emptyset$.  We claim that the intersection must be a point.  Let $l = \{t_l\} \times [s_0, s_1]$ and $m = \{t_m\} \times [s_2, s_3]$; without loss of generality, let $t_l \leq t_m$.  If $t_l = t_m$ and there is any intersection at all, this is actually a single line.  So we may also assume that $t_l < t_m$.  For all $r \in (s_1, s_3)$ and for all $t < t_m$, $(t, r)$ is in the interior of $\beta_{\text{left}}$.  The only possible with $l$ intersection is at an endpoint.  It cannot be both endpoints because this disconnects $l$.
    \par
    Take these horizontals and verticals as a starting point for the image of $f$.  Since there are finitely many horizontal and vertical segments, we may take $\cI$ to be the cardinality of the finite number of connected components.  By definition, $(3)$ is achieved.  All the horizontal segments of the boundary of $\beta_{\text{left}}$ are contained in $U_0$.  Therefore, there exists a tube about each compact horizontal segment contained in $U_0$.  If a horizontal segment $h$ is connecting two vertical segments $l$ and $m$ (perhaps one of these vertical segments is the boundary of $I \times I$, but this is no matter), there exists a line connecting a point in $l$ to a point in $m$ with nonzero slope that remains in $U_0$.  Tilting in the correct direction assures that the projection down onto the $s$ coordinate is injective. This achieves $(1)$ and $(4)$.
\end{proof}

\begin{prop}
    Let $A$ be the image of some function $f$ as in the Lemma 4.2.9.  There exists a strictly horizontal deformation retract of the leftward closure of $A$, $A_{\text{left}}$, onto $A$.  By strictly horizontal, we mean that, for all $\pi_s(\cD(t, s, r)) = s$ for all $(t, s, r) \in A_{\text{left}} \times I$.
\end{prop}
\begin{proof}
    The image of $f$ is homeomorphic to a finite coproduct of intervals.  $A_{\text{left}}$ has a connected component for each term in this coproduct.  It suffices to define the deformation retract on each connected component separately, and we may therefore assume that $A_{\text{left}}$ is connected and that $A$ is homeomorphic to an interval.
    \par
    To each horizontal strip at $s$ in  $A_{\text{left}}$, there is a unique point in $A$ on that horizontal strip, denote it by $(t_s, s)$.  The deformation retract acts on each horizontal slice by taking the interval $[0, t_s]$ and contracting it to the point $t_s$.  Parametrize $A$ as the image of a piecewise linear function not from the interval, but from $\pi_s(A)$.  Since this projection is a homeomorphism on the restriction to $A$, this is possible through a affine homeomorphism between $I$ and $\pi_s(A)$.  We can then take express $t_s$ as the result of some continuous functions on $A_{\text{left}}$, $t_s = \pi_t(f(\pi_s(s, t))$.
    \begin{equation}
        \cD(t, s, r) = (r(t_s - t) + t, s)
    \end{equation}
    With the expression of $t_s$ as a continuous function on $A$, this is continuous.  That it is the identity in the $s$ variable is manifest, and at $r = 1$, for all $(t, s)$ we have $\cD(t, s, 1) = (t_s, s)$, realizing the deformation retract.
\end{proof}

\begin{proof}
    \emph{(proof of $(1)$ in Theorem 4.2.7)}
    Construct $A$ as the image of a piecewise linear function as in Lemma 4.2.9.  $A_{\text{left}}$ has a finite number of connected components, so one may lift $\cH$ and construct homotopies one connected component at a time.  As such, we may assume that $A_{\text{left}}$ is connected .  Take the pullback of $\cD_\circ$ along $\cD_\bullet \circ \cH|$.  Fiber bundles are stable under pullback, so the pullback is a fiber bundle over $A_{\text{left}}$.  This space is connected, contractible and locally contractible, therefore this fiber bundle over $A_{\text{left}}$ is trivializable, that is $A_{\text{left}} \times_{X_p} U_p \setminus X_p$ is homeomorphic to $A_{\text{left}} \times L \times (0, 1)$.  We are also equipped with a map to the pullback from $A$ by the universal property of a pullback.  In a diagram, the situation is the following:
\begin{figure}[h]
    \[\begin{tikzcd}
    	&&& {A_{\text{left}} \times L \times (0, 1)} \\
    	\\
    	A &&& {A_{\text{left}} \times_{X_p} U_p \setminus X_p} &&& {U_p \setminus X_p} \\
    	\\
    	{A_{\text{left}}} &&& {A_{\text{left}}} &&& {X_p}
    	\arrow["\phi"{description}, from=1-4, to=3-4]
    	\arrow["{\phi^{-1}\omega}"{description}, from=3-1, to=1-4]
    	\arrow["\omega"{description}, from=3-1, to=3-4]
    	\arrow["\iota"{description}, hook, from=3-1, to=5-1]
    	\arrow["f"{description}, from=3-4, to=3-7]
    	\arrow["g"{description}, from=3-4, to=5-4]
    	\arrow[from=3-7, to=5-7]
    	\arrow[dashed, from=5-1, to=1-4]
    	\arrow["\cG"{description}, dotted, from=5-1, to=3-4]
    	\arrow["{\text{Id}}"{description}, from=5-1, to=5-4]
    	\arrow[from=5-4, to=5-7]
    \end{tikzcd}\]
    \caption{}
    \label{fig:myfigure}
\end{figure}
    One is in need of constructing $\cG: A_{\text{left}} \rightarrow A_{\text{left}} \times_{X_p} U_p \setminus X_p$ that is exiting with respect to the stratification by $P_{>p}$ given to the pullback by post-composition of the stratification of $U_p \setminus X_p$ with the natural map $f$.  $\cG$ must also make everything commute. Note that there are equalities $f \omega = \cH|$ and $g \phi = \mathsf{proj}_{A_{\text{left}}}$.  We will define the map $\cG$ by way of a map to the trivialization, which may be defined component wise.  For $g\cG$ to be the identity, it is clear that the would be map into the product $A_{\text{left}} \times L \times (0, 1)$ must be the identity in the $A_{\text{left}}$ coordinate.  In the $L \times (0, 1)$ coordinate, first retract onto $A$, then map into $A_{\text{left}} \times L \times (0, 1)$ by $\phi^{-1}\omega$, then project down to $L \times (0, 1)$.  This may be written explicitly as the composition
    $$
    \cG_0 := \bigg(\text{Id}, \mathsf{proj}_{L \times (0, 1)}\circ \phi^{-1}\circ \omega\circ \mathsf{ret}\bigg), \hspace{3em} \cG = \phi \circ \cG_0
    $$
    We now must check the commutativity requirements and the exiting requirements.  For commutativity, it is enough to check that $\cG_0 \circ \iota = \phi^{-1} \circ \omega$.  If this is so, we'd have
    $$
    \cG_0 \circ \iota = \phi^{-1} \circ \omega \iff \phi \circ \cG_0 \circ \iota = \omega \iff \cG \circ \iota = \omega.
    $$
    i.e. the top triangle of the left square commutes.  The bottom is guaranteed to commute already by our choice in the $A_{\text{left}}$ coordinate.
    $$
    \cG_0 \circ \iota = \bigg(\text{Id}, \mathsf{proj}_{L \times (0, 1)}\circ \phi^{-1}\circ \omega\circ \mathsf{ret}\bigg) \circ \iota = \bigg(\iota, \mathsf{proj}_{L \times (0, 1)}\circ \phi^{-1} \circ \omega \bigg)
    $$
    In the $L \times (0,1)$ coordinate, this is $\phi^{-1}\circ \omega$. In the second coordinate we have
    $$
    \iota = g \circ \omega = \mathsf{proj}_{A_{\text{left}}} \circ \phi^{-1} \circ \omega
    $$
    To see that $\cG$ is exiting it is enough to show that $\cG_0$ is.  This follows from the fact that the retract is strictly horizontal as the restriction of the deformation retract in the previous proposition.  Take $\tilde{\cH}| := f \circ \cG$ as the lift of $\cD_\bullet \circ \cH|$ into $U_p \setminus X_p$.  This is a lift of a \emph{restriction} of $\cH$.  To define the lift on all of $\cH$:
    \[
    \tilde{\cH} = 
    \begin{cases}
        \tilde{\cH}|(t, s) & \text{if } (t, s) \in A_{\text{left}} \\
        \cH(t, s) & \text{if } (t, s) \notin A_{\text{left}} 
    \end{cases}
    \]
    Since the outer rectangle in Figure 1 commutes, this agrees on $A$.  Since the retract of $A_{\text{left}}$ onto $A$ is strictly horizontal, this is exiting.
\end{proof}

\begin{proof}
    \emph{(proof of $(2)$ in Theorem 4.2.7)}  We continue to assume, without loss of generality, that $A_{\text{left}}$ is connected.  To produce a homotopy between $\cH, \tilde{\cH}$, it is sufficient to produce homotopies on the restrictions to $A_{\text{left}}$ rel $A$.  With this in mind, we simply refer to the restrictions of $\cH, \tilde{\cH}$ to $A_{\text{left}}$ by their names $\cH, \tilde{\cH}$ with no other indication.  To construct the homotopy, first extend $A_{\text{left}}^{ext}$ by attaching a rectangle to its left end.  Let $\epsilon < 0$ and sufficiently big, define
    $$
    A_{\text{left}}^{ext} = A_{\text{left}} \cup \bigg(\pi_s(A_{\text{left}}) \times [\epsilon, 0)\bigg)
    $$
    We can also extend the homotopies $\cH, \tilde{\cH}$ by defining $\cH_{ext}, \tilde{\cH}_{ext}$ on the newly adjoined rectangle by $\cH_{ext}(-, s) = \cH(0, s)$ for all $t \in  [\epsilon , 0)$.  Note that these extensions remain exiting along horizontal slices.  A homotopy between the extensions then restricts to a homotopy between $\cH$ and $\tilde{\cH}$.  The advantage of appending on this extra rectangle is that $A_{\text{left}}^{ext}$ is homeomorphic to $I \times I$ by a homeomorphism that exclusively stretches horizontally and vertically, and such that the image of $A$ is the right edge.  Because of this, the post-composition of $\cH_{ext}, \tilde{\cH}_{ext}$ with the homeomorphism remains exiting.  This reduces to the case that $A$ is the right edge and $A_{\text{left}}$ is the entire square.
    \par
    Define coordinates of $I \times I \times I$ as $(t, s, \alpha)$, with the $\alpha$ variable the homotopy variable, not taken to be $r$ to differentiate it from the variable we are using to retract $U_p$ onto $X_p$.  Both constructions of $\Phi$ and $\tilde{\Phi}$ are agnostic to which input $\cH$ or $\tilde{\cH}$, so we construct $\Phi$ first for $\cH$ then $\tilde{\Phi}$ is defined by replacing $\cH$ with $\tilde{\cH}$ in the formula.  After the construction is finished, we argue that the endoint homotopy $\cG$ is regardless of the input $\cH, \tilde{\cH}$.
    \par
    The deformation retract $\cD: U_p \times I \rightarrow U_p$ affords us one homotopy by post-composing $\cD$ with $\cH$.
    \[\begin{tikzcd}
    	{I \times I \times I} &&& {U_p \times I} &&& {U_p}
    	\arrow["{\cH \times \text{Id}}", from=1-1, to=1-4]
    	\arrow["\cD", from=1-4, to=1-7]
    \end{tikzcd}\]
    Let $\phi_1 := \cD \circ (\cH \times \text{Id})$.  The homotopy cube with axes labels looks like
    
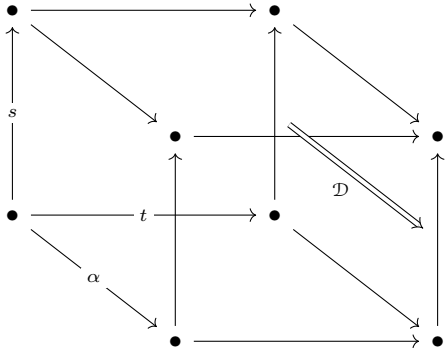
\[\begin{tikzcd}
    	\bullet &&& \bullet \\
    	\\
    	&& \bullet &&& \bullet \\
    	\bullet &&& \bullet \\
    	\\
    	&& \bullet &&& \bullet
    	\arrow[from=1-1, to=1-4]
    	\arrow[from=1-1, to=3-3]
    	\arrow[from=1-4, to=3-6]
    	\arrow[from=3-3, to=3-6]
    	\arrow["s"{description}, from=4-1, to=1-1]
    	\arrow["t"{description}, from=4-1, to=4-4]
    	\arrow["\alpha"{description}, from=4-1, to=6-3]
    	\arrow[""{name=0, anchor=center, inner sep=0}, from=4-4, to=1-4]
    	\arrow[from=4-4, to=6-6]
    	\arrow[from=6-3, to=3-3]
    	\arrow[from=6-3, to=6-6]
    	\arrow[""{name=1, anchor=center, inner sep=0}, from=6-6, to=3-6]
    	\arrow["\cD"', Rightarrow, shorten <=6pt, shorten >=6pt, from=0, to=1]
    \end{tikzcd}\]
    The back face is the original homotopy $\cH$.  This is exiting as long as $\cH$ is.  For all $s, t, \alpha \in (0,1)$, $\phi_1(s, t, \alpha)$ takes image in the same stratum as $\cH(s, t)$.  At $\alpha = 1$, the front face is trivially stratified; everything is mapped to the stratum $X_p$.  When we say that the input of the construction does not depend on $\cH$, we mean that we might put $\tilde{\cH}$ here and the homotopy at $\alpha = 1$ will be the same.  The right face takes the image of the right edge and deforms it to the cone locus.  Deform the square an epsilon amount to give us some space on the right side.
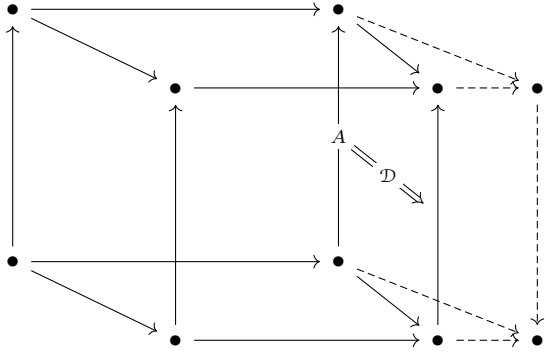
\[\begin{tikzcd}
	\bullet &&&& \bullet \\
	&& \bullet &&& \bullet & \bullet \\
	\\
	\\
	\bullet &&&& \bullet \\
	&& \bullet &&& \bullet & \bullet
	\arrow[from=1-1, to=1-5]
	\arrow[from=1-1, to=2-3]
	\arrow[from=1-5, to=2-6]
	\arrow[dashed, from=1-5, to=2-7]
	\arrow[from=2-3, to=2-6]
	\arrow[dashed, from=2-6, to=2-7]
	\arrow[dashed, from=2-7, to=6-7]
	\arrow[from=5-1, to=1-1]
	\arrow[from=5-1, to=5-5]
	\arrow[from=5-1, to=6-3]
	\arrow[""{name=0, anchor=center, inner sep=0}, "A"{description}, from=5-5, to=1-5]
	\arrow[from=5-5, to=6-6]
	\arrow[dashed, from=5-5, to=6-7]
	\arrow[from=6-3, to=2-3]
	\arrow[from=6-3, to=6-6]
	\arrow[""{name=1, anchor=center, inner sep=0}, from=6-6, to=2-6]
	\arrow[dashed, from=6-6, to=6-7]
	\arrow["\cD"{description}, Rightarrow, shorten <=6pt, shorten >=6pt, from=0, to=1]
\end{tikzcd}\]
    We intend to define the homotopy $\Phi$ on the closed subsets of $I \times I \times I$ given by the image of this homeomorphism and the closure of its complement in the cube.  The homotopy should be constant along the right face, so it makes sense to fill in the rest of the front face with the deformation of the image of $A$ onto its cone locus, as in the right face of the original homotopy cube.
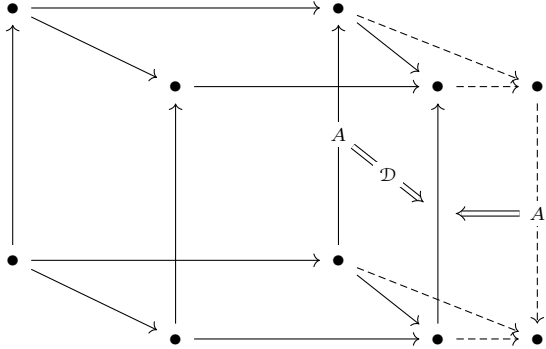
\[\begin{tikzcd}
	\bullet &&&& \bullet \\
	&& \bullet &&& \bullet & \bullet \\
	\\
	\\
	\bullet &&&& \bullet \\
	&& \bullet &&& \bullet & \bullet
	\arrow[from=1-1, to=1-5]
	\arrow[from=1-1, to=2-3]
	\arrow[from=1-5, to=2-6]
	\arrow[dashed, from=1-5, to=2-7]
	\arrow[from=2-3, to=2-6]
	\arrow[dashed, from=2-6, to=2-7]
	\arrow[""{name=0, anchor=center, inner sep=0}, "A"{description}, dashed, from=2-7, to=6-7]
	\arrow[from=5-1, to=1-1]
	\arrow[from=5-1, to=5-5]
	\arrow[from=5-1, to=6-3]
	\arrow[""{name=1, anchor=center, inner sep=0}, "A"{description}, from=5-5, to=1-5]
	\arrow[from=5-5, to=6-6]
	\arrow[dashed, from=5-5, to=6-7]
	\arrow[from=6-3, to=2-3]
	\arrow[from=6-3, to=6-6]
	\arrow[""{name=2, anchor=center, inner sep=0}, from=6-6, to=2-6]
	\arrow[dashed, from=6-6, to=6-7]
	\arrow[Rightarrow, shorten <=6pt, shorten >=6pt, from=0, to=2]
	\arrow["\cD"{description}, Rightarrow, shorten <=6pt, shorten >=6pt, from=1, to=2]
    \end{tikzcd}\]
    We now make this idea concrete.  Let $C_0$ denote the cube deformed by an epsilon amount.  That is, there is a continuous endomorphism of the cube
    \begin{equation}
        \omega(t, s, \alpha) = ((1-\alpha\epsilon)t, s, \alpha)
    \end{equation}
    that is a homeomorphism onto its image with inverse
    \begin{equation}
        \omega^{-1}(t, s, \alpha) = (\frac{1}{1-\alpha\epsilon}t, s, \alpha)
    \end{equation}
    
    $C_0$ is the image of this map.  Define $C_1$ as the closure of its complement in the cube.  These are given explicitly by
    \[
    C_0 = \bigg\{(t, s, \alpha) \in I \times I \times I : \alpha \leq \frac{-1}{\epsilon}(t-1) \bigg\}
    \]
    \[
    C_1 = \bigg\{(t, s, \alpha) \in I \times I \times I : \alpha \geq \frac{-1}{\epsilon}(t-1) \bigg\}
    \]
    We now define the homotopy $\Phi$ on each of these parts.
    \begin{equation}
        \Phi(t, s, \alpha) = 
        \begin{cases}
            \cD \circ (\cH \times \text{Id}) \circ \omega^{-1}(t, s, \alpha) & \text{if } (t, s, \alpha) \in C_0 \\
            \cD\big(\cH(1, s), \frac{-1}{\epsilon}(t-1)\big) & \text{if } (t, s, \alpha) \in C_1
        \end{cases}
    \end{equation}
    There are three things to check:
    \begin{enumerate}
        \item This is well defined.
        \item Replacing $\tilde{\cH}$ with $\cH$, this defines $\tilde{\Phi}$.  The restrictions of $\tilde{\Phi}$ and $\Phi$ to $\alpha = 1$ agree.
        \item For fixed $\alpha, s$, the resulting restriction in both $\Phi, \tilde{\Phi}$ is an exit path.
    \end{enumerate}
    Checking $(1)$ first, the boundary is when $\alpha = \frac{-1}{\epsilon}(t-1)$.
    \[
    \cD(\cH \times \text{Id}(\omega^{-1}(t, s, \frac{-1}{\epsilon}(t-1)))) = \cD(\cH \times \text{Id}\bigg(\big(\frac{1}{1-\frac{-1}{\epsilon}(t-1)\epsilon}\big)t, s, \frac{-1}{\epsilon}(t-1)\bigg) = 
    \]
    \[
    \cD(\cH \times \text{Id}(1, s, \frac{-1}{\epsilon}(t-1))) = \cD(\cH(1, s), \frac{-1}{\epsilon}(t-1))
    \]
    This is exactly the definition of $\Phi$ on $C_1$.  To see $(2)$, note that on $C_1$  $\Phi, \tilde{\Phi}$ agree everywhere because $\cH, \tilde{\cH}$ agree on the right edge of the square.  On $C_0$ at $\alpha = 1$, the images both $\Phi, \tilde{\Phi}$ are being deformed all the way to $X_p$, that is $\cD$ becomes the fiber bundle map $\cD|_{\setminus\frak{loc}}$ on $\tilde{\cH}$ and become $\cD|$ on $\cH$.  These agree by commutativity of the diagram in Theorem 4.2.7 used to construct $\tilde{\cH}$.
    \par
    Lastly, for $(3)$, fix some $s \in I$.  Since $\cH, \tilde{\cH}$ are exiting and $\omega$ is a homeomorphism that contracts in only the $t$ direction, $\Phi$ is certainly exiting on $C_0$.  On $C_1$, a fixed cross section looks like
    
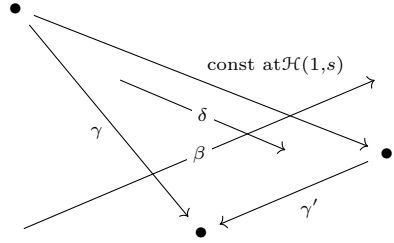
\[\begin{tikzcd}
    	\bullet \\
    	& {} & {} && {} \\
    	&&& {} & \bullet \\
    	{} && \bullet
    	\arrow["{\text{const at} \cH(1, s)}", from=1-1, to=3-5]
    	\arrow["\gamma"', from=1-1, to=4-3]
    	\arrow["\delta"{description}, from=2-2, to=3-4]
    	\arrow["{\gamma'}", from=3-5, to=4-3]
    	\arrow["\beta"{description}, from=4-1, to=2-5]
    \end{tikzcd}\]
    $\gamma, \gamma'$ are both deformations of $\cH(1, s)$ to the cone point.  Along $\delta$ it is constant at some point midway in the deformation.  We want to know if the trajectory along $\beta$ is exiting.  By proposition 4.1.12 however, all points in this cross section, except the bottom point when the deformation hits time 1, take image in the same stratum as $\cH(1, s)$.  Each cross section $\beta$ is constantly stratified except the one following the reverse of $\gamma'$, which exits immediately at the bottom point then goes on to be constantly stratified.
    \par
    For the homotopy with domain $A_{\text{left}}$, precompose with the inverse from $A_{\text{left}}^{ext}$ and restrict.  Since the inverse contracts only horizontally and vertically and the restriction of an exiting map remains exiting, this homotopy is exiting.
\end{proof}
\begin{proof}
    \emph{(proof of $(3)$ in Theorem 4.2.7)} Take $U_p$ small enough such that the fibers over all points in $I \times \{0\} \cap A_{\text{left}}$ are contained in $B$, as in Lemma 4.1.16.  This forces the image of $I \times \{0\} \cap A_{\text{left}}$ to be contained in $B$.  The homotopy $\Phi$ is along the fibers of the fiber bundle, so for all time $\alpha$, the bottom edge remains in $B$.
\end{proof}

\subsection{Replacing homotopies}

\begin{definition}
    A \emph{\textbf{chain}} (of length $n$) in a poset is a sequence of distinct elements $a_1 < \dots < a_n$.
\end{definition}
\begin{definition}
    For $P$ a finite poset, the \emph{\textbf{depth of an element}} $p$ is the length of the longest chain starting at a minimal element $m \in P$ and ending at $p$.  The \emph{\textbf{depth of a poset}} is the value $\mathsf{max}_{p \in P}\{\mathsf{depth}(p)\}$.  Equivalently, one can take this to be the length of the longest chain in $P$.
\end{definition}
\begin{definition}
    The \emph{\textbf{chain of length $n$}}, is the poset $1 < 2 < \dots < n$ and denoted $[n]$
\end{definition}
\begin{prop}
    For $P$ a finite poset, for some $n \geq 0$ there exists a continuous map $d: P \rightarrow [n]$ such that $p < q \implies d(p) < d(q)$, where $n = \mathsf{depth}(P)$.  If $\cS: X \rightarrow P$ is a conically smooth stratification, then so is $f \circ \cS$.  The strata of $X$ given by the poset $[n]$ are $X_k = \coprod_{p \in P : \mathsf{depth}(p) = k} X_p$.
\end{prop}
\begin{proof}
    Induct on the depth of the poset.  For (connected) posets with depth one, this is just a singleton.  We may take $[n] = [1]$ and the map $d$ is the identity.
    \par
    Suppose that this holds for all posets of depth less than $n$.  Fix a poset $P$ of depth $n$ and consider the new poset $P_0$ acquired by removing all the maximal elements form $P$.  This poset has depth $n-1$ and is equipped by the inductive hypothesis with a continuous map $d_0: P_0 \rightarrow [n-1]$ such that $p < q$ implies that $d_0(p) < d_0(q)$.  $d_0$ extends to a map $d: P \rightarrow [n]$ by defining for all maximal elements $m \in P$, $d(m) = n$.
\end{proof}
\begin{remark}
    The utility of this observation is that, given a conically smooth stratification $\cS: X \rightarrow P$, we may post-compose $\cS$ with the map $d$ and assume that $X$ is stratified by a linearly ordered poset.
\end{remark}

\begin{prop}
    Let $\cH$ be an exit path homotopy rel endpoints in a stratified space $X$ whose stratification is conically smooth.  It is homotopic to an immediately exiting homotopy.
\end{prop}
\begin{proof}
    Without loss of generality, assume that the stratifying poset of $X$ is $[n]$ for some $n$ and first construct a homotopy $\Phi: \cH \rightarrow \cH_0$ where $\cH_0$ takes image entirely in a single stratum.  This is done by repeated applications of Theorem 4.2.7.  One 'peels' the homotopy off of each stratum, one at a time.  Note that while this homotopy does not preserve the left edge, we may take it to preserve the right edge by taking $U_p$ to be small enough at every step.  Explicitly, if $\cH(1, s) = y$ for all $s \in I$, and $y \in U_p$, since $y$ is in the top most strata and we have not finished the process yet, it must be so that $y \in U_p \setminus X_p$. .  Deform $U_p$ until $y \notin U_p$.  The resulting homotopy achieved by concatenating all the peels satisfies $\tilde{\cH}(1, s) = y$ for all $s \in I$.
    \par
    Let $\cH(0, s) = x$ for all $s \in I$.  Part $(3)$ of Theorem 4.2.7 guarantees that when we peel, we do not peel things too far away $x$; we can take $\Phi$ so that for every $s, \alpha$, $\Phi(0, s, \alpha) \subseteq B$ for $B$ an arbitrary open set about $x$. Take such an open set as a basic of the form $V \times C(L)$, where $V$ is a contractible vector space.  There is a straight line deformation retract of $C(L)$ onto its cone point by $\cD: C(L) \times I \rightarrow C(L)$, $\cD([l, \alpha], t) = [l, \alpha t]$.  Combine this with the deformation retract of the vector space $V$ to the origin yields a deformation retract of the open basic $V \times C(L)$ onto $x$.  The resulting path from $x$ to each point along this deformation is immediately exiting.  In a picture:
    
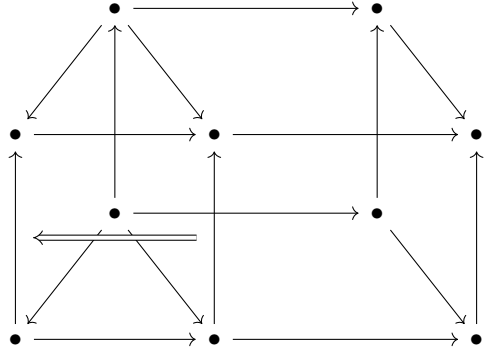
\[\begin{tikzcd}
    	& \bullet &&& \bullet \\
    	\\
    	\bullet && \bullet &&& \bullet \\
    	& \bullet &&& \bullet \\
    	\\
    	\bullet && \bullet &&& \bullet
    	\arrow[from=1-2, to=1-5]
    	\arrow[from=1-2, to=3-1]
    	\arrow[from=1-2, to=3-3]
    	\arrow[from=1-5, to=3-6]
    	\arrow[from=3-1, to=3-3]
    	\arrow[from=3-3, to=3-6]
    	\arrow[from=4-2, to=1-2]
    	\arrow[from=4-2, to=4-5]
    	\arrow[from=4-2, to=6-1]
    	\arrow[from=4-2, to=6-3]
    	\arrow[from=4-5, to=1-5]
    	\arrow[from=4-5, to=6-6]
    	\arrow[""{name=0, anchor=center, inner sep=0}, from=6-1, to=3-1]
    	\arrow[from=6-1, to=6-3]
    	\arrow[""{name=1, anchor=center, inner sep=0}, from=6-3, to=3-3]
    	\arrow[from=6-3, to=6-6]
    	\arrow[from=6-6, to=3-6]
    	\arrow[Rightarrow, shorten <=6pt, shorten >=6pt, from=1, to=0]
        \end{tikzcd}\]
    Since each deformation is immediately exiting, this remains immediately exiting for all fixed $\alpha$, $s$.  The new front face of the homotopy cube is an immediately exiting homotopy.
\end{proof}
\begin{cor}
    Any exit path $f$ is homotopic to an immediately exiting path $f'$.
\end{cor}

We've shown the special case of lifting homotopies when they are immediately exiting in Proposition 4.0.3.  This isn't quite enough to show that $\Pi$ is a right fibration.  One needs to show that any lift of a path is well defined.  That is, if $f, g$ are homotopic in the quotient $X/G$, the lifts obtained after lifting $f(1), g(1)$ are homotopic in $X$.  The above construction takes $f \simeq_{\text{exit}} g$ and constructs homotopies 
$$
f \simeq_{\text{exit}} f' \simeq_{\text{im. exit}} g' \simeq_{\text{exit}} g
$$
The immediately exiting homotopy witnessing $f'$ homotopic to $g'$ lifts to a homotopy between $\tilde{f}'$ and $\tilde{g}'$, but one also wants to know that the homotopy between $f, f'$ also lifts to a homotopy between $\tilde{f}$ and $\tilde{f}'$.  Fortunately, the motivating case in proposition 4.0.3 makes an assumption that is far too strong.  Generally speaking, one \emph{expects} to be able to lift homotopies by way of covering space theory.  That is, should one write down any reasonable homotopy whatsoever, the boundary of its strata will be well-behaved enough to carry on.  Immediately exiting homotopies are exceedingly reasonable.  So is the constructed homotopy between $f$ and it's immediately exiting counterpart $f'$.  It is by way of the deformation retract of $U_p$ onto the stratum $X_p$ at every stage.  
\par
\begin{lemma}
    Let $\cH: I \times I \rightarrow X/G$ be an exit path homotopy.  The induced stratification of $I \times I$ has three strata given by $p < q < r$ in the stratifying poset $P_G^{\op}$.
    \begin{enumerate}
        \item The $p$ stratum is $I \times \{0\}$
        \item The $q$ stratum is $\{0\} \times (0, 1]$
        \item The $r$ stratum is $(0, 1] \times (0, 1]$
    \end{enumerate}
    A choice of lift for $(1, 1)$ determines a lift $\tilde{\cH}: I \times I \rightarrow X$.
\end{lemma}
\begin{proof}
    Lift the first stratum first.  $(0, 1] \times (0, 1]$ takes image in a single stratum and is contractible.  Therefore, a choice of lift for $(1, 1)$ determines a map $\tilde{\cH}_{(0, 1] \times (0, 1]}: (0, 1] \times (0, 1] \rightarrow X_r$ that lifts $\cH_r$.
    \par
    This establishes partial lifts for all the horizontal and vertical cross sections of the square.  Indeed, $\cH_{I \times \{s\}}$, $\cH_{\{t\} \times I}$ are continuous paths into $X/G$, and the lift $\tilde{cH}_r$ fixes a choice of the endpoint of all of these except $s = 0, t = 0$.  After fixing the endpoint, since all these paths are exiting, this determines lifts of the entire path.  Furthermore, by uniqueness of lifts, $\tilde{\cH}_r|_{\{t\} \times (0, 1]}$ agrees with the lift of $\cH|_{\{t\} \times I}$ on $(0, 1]$ in the second coordinate.  We may argue the same logic in the $s$ direction in a way that it entirely symmetric.  Each of these horizontal and vertical extensions of $\tilde{\cH}_r$ the lift on the $r$ stratum defines a set map $\tilde{\cH}_{\setminus (0, 0)}$ from the square remove the origin to $X$. We now endeavor to show that this set map is continuous.  To do this, it suffices to show that it is continuous at some point $(t, 0)$ with $t \neq 0$.  It is continuous at $(t, s)$ for all $t \neq 0, s\neq 0$ by continuity of $\tilde{\cH}_r$, and the case of $(0, s)$, $s \neq 0$ is entirely symmetric.
    \par
    Let $x := \tilde{\cH}_{\setminus (0, 0)}(t, 0)$ and let $U$ be some arbitrary open subset about $x$.  As in lemma 3.1.4, find some $V \subseteq U$ satisfying the properties of inclusion, disjointness and symmetry.  $\cH(t, 0) = [x]$ in $X/G$ and $\pi(V)$ is an open set about $[x]$.  Continuity of $\cH$ implies that there exists some connected open box $B$ about $(t, 0)$ such that $\cH(B) \subseteq \pi(V)$, or rather an open box intersected with the square.  It is convenient to be specific in what $B \cap I \times I$ is, take it to be of the form $(t - \epsilon, t + \epsilon) \times [0, \epsilon) \subseteq I \times I$.  This is open and connected, and the removal of the bottom edge remains open and connected.  Define 
    $$
    B_{\setminus} = (t - \epsilon, t + \epsilon) \times [0, \epsilon) \subseteq I \times I \setminus I \times \{0\} = (t - \epsilon, t + \epsilon) \times (0, \epsilon)
    $$
    This is an open set in $(0, 1] \times (0, 1]$.  Since $\tilde{\cH}_{\setminus (0, 0)}$ is continuous on this subset, its image is a connected set in $\pi^{-1}(\pi(V)) = \coprod g \cdot V$, which means that it must take image in one cofactor of the coproduct.  We claim that the cofactor it takes image in is $V$ itself.  $\tilde{\cH}_{\setminus (0, 0)}$ is continuous on the vertical strip at $t$, denote this as $f_t = \tilde{\cH}_{\setminus (0, 0)}|_{\{t\} \times I}$.
    Since $x = f_t(0) \in V$, $f_t((0, \epsilon)) \subseteq g \cdot V$ for $g$ not in the stabilizer of $x$, produces a contradiction to continuity of $f_t$.  Since $f_{t_0}((0, \epsilon)) \subseteq V$ also holds and $f_{t_0}$ is also continuous for all $t_0$ in the $t$ coordinate of $B \cap I \times I$, it must be so that $f_{t_0}(0) \in V$
    \par
    What remains is to fill in the last point, $(0, 0)$.  If $\cH(0, 0) = [x]$, let $\pi(V)$ by some symmetric, disjoint set about a lift of $[x]$.  $\cH$ is continuous, this implies there is an open box
    $$
    B = [0, \epsilon) \times [0, \epsilon) \subseteq I \times I
    $$
    Such that $\cH(B) \subseteq \pi(V)$.  Now, $B \setminus \{(0, 0)\}$ is a connected open subspace of $I \times I \setminus \{(0, 0)\}$.  Since we've shown that our lift is already continuous on this subspace, it takes image in a single cofactor of $\pi^{-1}(\pi(V)) \coprod g \cdot V$.  Without loss of generality, we may assume that it takes image in $V$.  Define $\tilde{\cH}(0, 0) = x$ for $x$ the unique lift of $[x]$ in $V$.  
    \par
    We finish the proof by showing continuity at this point. Let $U$ be some open subset about $x$.  There exists $W \subseteq U \cap V$ that is symmetric, disjoint in its orbit.  By continuity of $\cH$, there exists connected, open $W_0$ about $(0, 0)$ such that $\cH(W_0) \subseteq \pi(W)$.  We may also take $W_0 \subseteq B$ by taking an intersection if needed and taking a connected component.  Again, the restriction to $I \times I \setminus \{(0, 0)\}$ is continuous.  This implies $W_0$ is mapped to a single cofactor of $\coprod g \cdot W$.  But $W_0 \subseteq B$ implies that its image is a subset of the image of $B$.  The image of $B$ is contained in $V$, forcing the image of $W_0$ to be contained in $W \subseteq V$.
\end{proof}
\begin{observation}
    Nowhere in this proof do we use the fact that the strata $p, q$ are different.  This proof also goes of $p = q$ and have stratum $I \times \{0\} \cup \{0\} \times (0, 1]$ in $I \times I$.  We also might reorient the edges of the square and pick a different starting lift, or even consider a homeomorphism from the square to a rectangle with a similar stratification!
\end{observation}
\begin{theorem}
    Let $\cH$ be an exit path homotopy between exit paths $f, g$.  Let $\Phi$ be a homotopy from $\cH$ to $\tilde{\cH}$ as in Proposition 4.3.6 where $\tilde{\cH}$ is immediately exiting.  Let $\cF, \cG$ denote the homotopy relating $f \simeq_{\text{exit}} f'$ and $g \simeq_{\text{exit}} g'$ respectively.  Then a choice of lift for $f(1) = f'(1)$ determines a lift of $\cF$.  Likewise for $g(1) = g'(1)$ and $\cG$.
\end{theorem}
\begin{proof}
    It suffices to lift only $\cF$.  The homotopy $\cF$ is a concatenation of homotopies between intermediary paths.  A lift of $\cF$ follows from lifts of each of these individual homotopies, all of which take the same form: they are by way of the deformation retract, modulo a deformation on each end of the cube.  Up to some deformation of the cube, this looks like
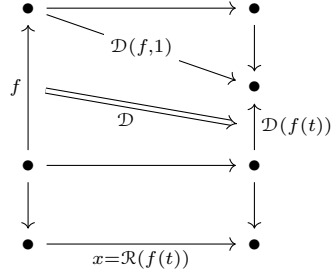
\[\begin{tikzcd}
	\bullet &&& \bullet \\
	&&& \bullet \\
	\bullet &&& \bullet \\
	\bullet &&& \bullet
	\arrow[from=1-1, to=1-4]
	\arrow["{\cD(f, 1)}"{description}, from=1-1, to=2-4]
	\arrow[from=1-4, to=2-4]
	\arrow[""{name=0, anchor=center, inner sep=0}, "f", from=3-1, to=1-1]
	\arrow[from=3-1, to=3-4]
	\arrow[from=3-1, to=4-1]
	\arrow[""{name=1, anchor=center, inner sep=0}, "{\cD(f(t))}"', from=3-4, to=2-4]
	\arrow[from=3-4, to=4-4]
	\arrow["{x = \cR(f(t))}"', from=4-1, to=4-4]
	\arrow["\cD"',Rightarrow,shorten <=6pt, shorten >=6pt, from=0, to=1]
\end{tikzcd}\]
    We provide this picture as a heuristic.  A typical homotopy in the concatenation is ultimately defined on 3 different parts of the square.  All three parts meet at boundaries homeomorphic to intervals, and the image of all the boundaries are \emph{exit paths}. Indeed, the shared bottom edge is deformation of $f(0)$ to the cone locus, this is an immediately exiting path ending at $f(0)$.  The diagonal shared top edge is the deformation of $f(1)$ to the cone locus, also an immediately exiting path.  It suffices then to lift all parts separately. Since path lifting is unique given a choice of lift of the endpoint, the corresponding lifts will agree on the boundaries and be amenable to any kind of gluing.  A choice of $f(1)$ is to lift the top triangle and the middle deformed square.  The lift of the middle deformed square determines a choice of a lift of $f(0)$, which should determine a lift of the bottom rectangle, the deformation in $V \times C(L)$.
    \par
    Considering the stratification of the top triangle, it is constantly stratified by the stratum in which $f(1)$ lies everywhere except the bottom right point, which takes image in the cone locus.  This is an immediately exiting triangle so to speak and lifts by the fact that immediately exiting homotopies lift.  The stratification of the bottom rectangle is exactly of the form in Lemma 4.3.8, up to a shift in what edges are the lower dimensional stratum and what the choice of initial lift is.  By this lemma it too lifts with a choice of lift for $f(0)$.
    \par
    For these purposes, in considering the deformed middle square, it suffices to work with the undeformed square:
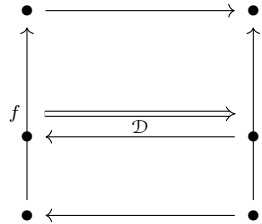
\[\begin{tikzcd}
	\bullet &&& \bullet \\
	\\
	\bullet &&& \bullet \\
	\bullet &&& \bullet
	\arrow[from=1-1, to=1-4]
	\arrow[from=3-4, to=3-1]
	\arrow[""{name=0, anchor=center, inner sep=0}, "f", from=4-1, to=1-1]
	\arrow[""{name=1, anchor=center, inner sep=0}, from=4-4, to=1-4]
	\arrow[from=4-4, to=4-1]
	\arrow["\cD"', Rightarrow, shorten <=6pt, shorten >=6pt, from=0, to=1]
\end{tikzcd}\]
    All the paths from the right edge to the left edge are either constantly stratified or immediately exiting!  In either case, the strata take the same form as in Lemma 4.3.8, up to a homeomorphism to a rectangle.  We may then lift strating with the top rectangle and proceed down one rectangle at a time.
\end{proof}

This completes what is perhaps the main moral achievement paper.  At the beginning of these proceedings, we commented that the main obstruction to lifting homotopies using covering space theory and lifting stratum by stratum was that there was no analog to proposition 2.1.13 in the 2-dimensional setting of $I \times I$.  The induced strata in $I \times I$ need not have any particularly nice form, as is the case for the interval.  The preceding proposition says exactly that, given two homotopic paths $f, g$, one can construct a new homotopy between them with strata that \emph{are} of a nice form.  Proposition 4.3.9 is as close to an analog of proposition 2.1.13 in the two dimensional as one might hope.


\section{The Main Theorem}
\subsection{$\mathsf{Enter}(M)$ as a pullback}
\begin{definition}
    A functor $p: E \rightarrow B$ is a \emph{\textbf{right fibration}} if for each $e \in E$ and arrow $f: b \rightarrow pe$ in $B$, there exists a unique lift $g: e' \rightarrow e$ in $E$.
\end{definition}
Some areas of the literature call this a \emph{discrete fibration}.  We choose this notation to evoke the corresponding dual, a left fibration.  The terminology 'discrete' is in particular used in Riehl and Loregian \cite{Riehl}, in which the following theorem may be found.
\begin{theorem}
    Let $\mathsf{DFib}(B)$ denote the category of discrete fibrations over $B$, defined to be the full subcategory of the comma category $\Cat/B$.  There is an equivalence of categories $\mathsf{DFib}(B) \simeq \fun(B^{\op}, \Sets)$.
\end{theorem}
\begin{definition}
    The \emph{\textbf{enter path category}}, $\mathsf{Enter}(X)$, associated to a stratification $\cS:X \rightarrow P$ is the category $\exit(X)^{\op}$.
\end{definition}
One may also write out explicitly what an enter path might be: $f: I \rightarrow X$ such that for all $t \leq s$, $\cS(f(t)) \geq \cS(f(s))$; one flips the direction of the inequality.  This category has as its objects the points of $X$ and as its morphisms these so called enter paths up to enter path homotopy.  To indicate the intention of moving towards proving the main theorem, we change notation from a space $X$ to a manifold $M$.
\begin{theorem}
    If $G$ is a finite group acting smoothly on a manifold $M$, then $\Pi: \exit(M) \rightarrow \exit(M/G)$ is a right fibration.
\end{theorem}
\begin{proof}
    Let $[x], [y]$ be two objects of $\exit(M/G)$ and $f$ a path class with $f(0) = [x], f(1) = [y]$.  A choice $y_0$ for $[y]$ determines a lift $\tilde{f}$ of $f$.  To see that this choice $\tilde{f}$ is well defined, suppose that $f \simeq_{\text{exit}} g$, $g$ is another representative of the path class of $f$. Replace the homotopy between $f, g$ with homotopies from $f$ to $f'$, the immediately exiting homotopy from $f'$ to $g'$, and the homotopy from $g'$ to $g$.  By Proposition 4.3.10 and Proposition 4.0.3, all of these lift to witness a homotopy from $\tilde{f}$ to $\tilde{g}$ in $M$.
\end{proof}

The right fibration $\Pi$ induces a functor $\Omega: \Enter(M/G) \rightarrow \set$ that operates as follows 
\[\begin{tikzcd}
	{[x]} &&& {\Omega[x] = \{\text{lifts of } [x]\}} \\
	\\
	{[y]} &&& {\Omega[y] = \{\text{lifts of } [y]\}}
	\arrow[curve={height=-18pt}, from=1-1, to=1-4]
	\arrow["f"', from=1-1, to=3-1]
	\arrow["{\Omega f}", from=1-4, to=3-4]
	\arrow[curve={height=-18pt}, from=3-1, to=3-4]
\end{tikzcd}\]
$f$ is an \emph{enter} path and corresponds to an \emph{exit} path from $[y]$ to $[x]$. For some lift $x_1$ of $[x]$, $\exists ! y_1$ and $f_0$ an exit path from $y_1$ to $x_1$. $\Omega f (x_1) = y_1$.
We can also consider $\Pi^{\op}: \Enter(M) \rightarrow \Enter(M/G)$ and consider a functor $\Omega_*: \Enter(M) \rightarrow \set_*$ making the square
\[\begin{tikzcd}
	{\Enter(M)} &&& {\mathsf{Sets}_*} \\
	\\
	{\Enter(M/G)} &&& {\mathsf{Sets}}
	\arrow["{\Omega_*}", from=1-1, to=1-4]
	\arrow["{\Pi^{\op}}"', from=1-1, to=3-1]
	\arrow["forget", from=1-4, to=3-4]
	\arrow["\Omega"', from=3-1, to=3-4]
\end{tikzcd}\]
commute.  The functor $\Omega_*$ is given below.
\[\begin{tikzcd}
	x &&& {\Omega_*x = \{\text{Orbit}(x), x\}} \\
	\\
	y &&& {\Omega_*y = \{\text{Orbit}(y), y\}}
	\arrow[curve={height=-18pt}, from=1-1, to=1-4]
	\arrow["f"', from=1-1, to=3-1]
	\arrow["{\Omega_*f}", from=1-4, to=3-4]
	\arrow[curve={height=-18pt}, from=3-1, to=3-4]
\end{tikzcd}\]
For the function $\Omega_*f$, the enter path $f$ determines one assignment, $\Omega_* f(x) = y$.  For the rest, note that $\text{Orbit}(x)$ and $\text{Orbit}(y)$, are transitive $\gset$, not just sets.  In particular, one should construct the functor $\Omega_*$ so that the corresponding morphism $\Omega_*f$ is $G$-equivariant.  Any $x' \in \text{Orbit}(x) = g \cdot x$ for some $g \in G$.  Define $\Omega_* f (x') = g \cdot y$.  One needs this definition to be well defined, i.e. if $x '  = g \cdot x = h \cdot x$, then $g \cdot y = h \cdot y$.  These both hold if and only if $h^{-1}g \in G_x, h^{-1}g \in G_y$.  If $x, y$ have an enter path between them, it is necessarily so that $[G_x] \subseteq [G_y]$.  $h^{-1}g \in G_x$ is provided by assumption and one wants to show that it is an element of $G_y$.  Therefore, we are in a position that requires us to show that $G_x \subseteq G_y$.  The next proposition provides this.
\begin{prop}
    Let $f: I \rightarrow M$ be an enter path with respect to the stratification of $M$ by $P_G^{\op}$.  It is also an enter path with respect to the stratification $S_G^{op}$
\end{prop}
\begin{proof}
    Let $f(0) = x$ and $f(1) = y$ with nonequal stabilizers $G_x$ and $G_y$.  We may assume that the stratum associated to $[G_x]$ and to $[G_y]$ are the only strata in which $f$ takes image and show that $G_x \subseteq G_y$.  Indeed, should this be true, in the general case one may look at the stratification of $I$ of the form $[0, a_1], (a_1, a_2], \dots, (a_{n-1}, 1]$ as in Proposition 2.1.13 and break it up into segments homeomorphic to the interval.  Apply this result first to $[a_1, a_2]$, then to $[a_2, a_3]$ etcetera.
    \par
    Since $f$ is an enter path with respect to the stratification of $M$ by $P_G^{\op}$, we know that $[G_x] \subseteq [G_y]$.  It follows that $G_x$ is conjugate to a subgroup of $G_y$.  Again invoking Proposition 2.1.13, suppose that the stratification of the interval is of the form $[0, \alpha], (\alpha, 1]$.  It also follows that for all $t \in [0, \alpha]$, $G_{f(t)}$ is conjugate to $G_x$.  Consider $f|_{[0, \alpha]}$.  Post composing this by the stratification $M$ by $S_G^{op}$, restrict the codomain of this continuous function to those elements of the poset that are conjugate to $G_x$.  This restricted codomain is completely disconnected and implies that the image of $f|$ in this codomain is constant.  Applying the same argument to $(\alpha, 1]$, we then have that for all $t \in [0, \alpha]$, $G_{f(t)} = G_x$ and for all $s \in (\alpha, 1]$, $G_{f(s)} = G_y$.
    \par
    Now, for the sake fo contradiction, suppose that $G_x$ is \emph{not} a subgroup of $G_y$.  This means that in the poset $S_G^{\op}$, there is not a relation $G_y \rightarrow G_x$.  Restrict the codomain of all of $f$.  By the previous argument, this is a restriction to two points, $G_x$ and $G_y$.  If there is no relation $G_y \rightarrow G_x$, this restricted codomain is again totally disconnected, witnessing a contradiction.
\end{proof}
\begin{cor}
    $\Omega_*$ is a functor.  $\Omega_*$ factors through $\cO_{G, *}$.
\end{cor}

\begin{prop}
The diagram
\[\begin{tikzcd}
	{Enter(M)} &&& {\mathsf{Sets}_*} \\
	\\
	{Enter(M/G)} &&& {\mathsf{Sets}}
	\arrow["{\Omega_*}", from=1-1, to=1-4]
	\arrow["{\Pi^{\op}}"', from=1-1, to=3-1]
	\arrow["\mathsf{forget}", from=1-4, to=3-4]
	\arrow["\Omega"', from=3-1, to=3-4]
\end{tikzcd}\]
is a pullback
\end{prop}
\begin{proof}
    Firstly, see that this commutes.  An object $x \in \mathsf{Enter}(M)$ is mapped to the pointed set $(\text{Orbit}(x), x)$,  the forgetful functor takes it to $\text{Orbit}(x)$.  $\Pi^{\op}$ of $x$ maps it to its quotient $[x]$, which is taken to $\{\text{lifts of } [x]\} = \text{Orbit}(x)$.  This is commutativity on the level of objects.
    \par
    On the level of morphisms, along the top path an enter path $f$ starting at $x$ and ending at $y$ results in a set function $\mathsf{forget}\Omega_* f: \text{Orbit}(x) \rightarrow \text{Orbit}(y)$ by $g \cdot x \rightarrow g \cdot y$.  Along the bottom path, $f$ is mapped to the enter path $\pi \circ f$ starting at $[x]$ and ending at $[y]$.  $\Omega(\pi \circ f)$ is a set function from $\text{Orbit}(x)$ to $\Orbit(y)$ that acts as follows.  The reverse of $\pi \circ f$ is an exit path ending at $[x]$.  Fix any lift $x_0$, this determines an exit path that lifts $[y]$ to some $y_0$.  In particular, fix $x$ as a lift for $x$.  Since $\pi \circ f$ is lifted by $f$, the corresponding lift of $[y]$ is $y$ itself.  Fix any $g \cdot x$ then as a lift of $[x]$, then the reverse of $g \cdot f$ with endpoint $g \cdot y$ serves as the unique exit path with lift starting at $g \cdot x$.  This is the same set map as $\mathsf{forget}\Omega_* f$.
    \par
    The argument assuring commutativity is useful in computing the pullback.  In doing so, we will see that it is an exact description of $\mathsf{Enter}(M)$ and finish the proof.  Objects consist of pairs $([x], (A, a))$ such that $\Omega[x] = A$.  This holds if $A = \text{Orbit}(x)$ and $a = x_0$ for some $x_0 \in \text{Orbit}(x)$.  This is just a point in $M$, along with its orbit data.  Clearly, there is a one to one correspondence with such objects and points of $M$.  The objects of the pullback are then exactly the objects of $\Enter(M)$.
    \par
    As for morphisms, a morphism in the pullback is a pair of morphisms that map to the same morphism in $\set$. Specifically, it is an enter path $f$ starting at $[x]$ and ending at $[y]$ along with the data of a map of pointed sets $g: (\text{Orbit}(x), x_1) \rightarrow (\text{Orbit}, y_1)$ such that $\mathsf{forget} g$ is the same as $\Omega f$. In particular, $g$ is $G$-equivariant.  There is a one-to-one correspondence between lifts $f$ and lifts of $[x]$ by the uniqueness of lifting exit paths.  By the characterization of commutativity, if $\mathsf{forget} g = \Omega f$, then the endpoint for the lift of $f$ at $x_1$ must be $y_1$.  This describes the morphisms of the pullback as exactly the enter paths in $\mathsf{Enter}(M)$.
\end{proof}

\begin{theorem}
    \[\begin{tikzcd}
	{\Enter(M)} && {\mathcal{O}_{G, *}} \\
	\\
	{\Enter(M/G)} && {\mathcal{O}_G}
	\arrow[from=1-1, to=1-3]
	\arrow[from=1-1, to=3-1]
	\arrow[from=1-3, to=3-3]
	\arrow[from=3-1, to=3-3]
    \end{tikzcd}\]
    is a pullback.
\end{theorem}
\begin{proof}
    $\Omega$ and $\Omega_*$ both factor through $\cO_G, \cO_{G,*}$ and the forgetful functors respectively.  That is, we have a diagram
\[\begin{tikzcd}
	&& {\mathcal{O}_{G, *}} \\
	{\Enter(M)} \\
	&&& {\set_*} \\
	&& {\mathcal{O}_G} \\
	{\Enter(M/G)} \\
	&&& {\set}
	\arrow[from=1-3, to=3-4]
	\arrow[from=1-3, to=4-3]
	\arrow[from=2-1, to=1-3]
	\arrow[from=2-1, to=3-4]
	\arrow[from=2-1, to=5-1]
	\arrow[from=3-4, to=6-4]
	\arrow[from=4-3, to=6-4]
	\arrow[from=5-1, to=4-3]
	\arrow[from=5-1, to=6-4]
\end{tikzcd}\]
    The bottom square in the triangle is a pullback by Proposition 5.1.7.  The right square is a pullback by Proposition 1.0.11.  By the two of three lemma for pullbacks, the square
    \[\begin{tikzcd}
	{\Enter(M)} && {\mathcal{O}_{G, *}} \\
	\\
	{\Enter(M/G)} && {\mathcal{O}_G}
	\arrow[from=1-1, to=1-3]
	\arrow[from=1-1, to=3-1]
	\arrow[from=1-3, to=3-3]
	\arrow[from=3-1, to=3-3]
\end{tikzcd}\]
is a pullback. 
\end{proof}

\printbibliography

\end{document}